\documentclass[aop,MSNbibl,seceqn,dvips]{arximspdf}

\doi{10.1214/12-AOP745} 
\volume{41}
\issue{4}
\pubyear{2013}
\firstpage{2791}
\lastpage{2819}

\makeatletter
\newcommand{\rrvert}{\vert}
\newcommand{\llvert}{\vert}

\newtheorem{theorem}{Theorem}[section]

\newtheorem{lemma}[theorem]{Lemma}

\newproclaim{remark}[theorem]{Remark}
\newproclaim{example}[theorem]{Example}

\newcommand{\prt}{\partial}
\newcommand{\eps}{\varepsilon}
\newcommand{\D}{\mathbb{D}}
\renewcommand{\P}{\mathbb{P}}
\newcommand{\E}{\mathbb{E}}
\newcommand{\R}{\mathbb{R}}
\newcommand{\RR}{\mathbb{R}}
\newcommand{\QQ}{{\mathbb Q}}
\newcommand{\Z}{\mathbb{Z}}
\newcommand{\ZZ}{\mathbb{Z}}
\newcommand{\A}{{\mathcal A}}
\newcommand{\B}{{\mathcal B}}
\newcommand{\M}{{\mathcal M}}

\newcommand{\CC}{{\mathcal C}}
\newcommand{\EE}{{\mathcal E}}
\newcommand{\FF}{{\mathcal F}}
\newcommand{\GG}{{\mathcal G}}

\newcommand{\vphi}{\varphi}
\newcommand{\bone}{{\mathbf1}}
\newcommand{\dist}{\operatorname{dist}}

\newcommand{\df}{\stackrel{\mathrm{df}}{=}}

\makeatother

\begin{document}
\begin{frontmatter}

\title{Reflecting random walk in fractal domains\thanksref{T1}}
\runtitle{Reflecting random walk in fractal domains}

\thankstext{T1}{Supported in part by NSF Grants DMS-09-06743
and DMR-1035196, and by Grant N N201 397137, MNiSW, Poland.}

\begin{aug}
\author[A]{\fnms{Krzysztof} \snm{Burdzy}\ead[label=e1]{burdzy@math.washington.edu}}
\and
\author[A]{\fnms{Zhen-Qing} \snm{Chen}\corref{}\ead[label=e2]{zqchen@uw.edu}}
\runauthor{K. Burdzy and Z.-Q. Chen}
\affiliation{University of Washington}
\address[A]{Department of Mathematics\\
University of Washington\\
Box 354350\\
Seattle, Washington 98195\\
USA\\
\printead{e1}\\
\hphantom{E-mail: }\printead*{e2}} 
\end{aug}

\received{\smonth{5} \syear{2011}}
\revised{\smonth{2} \syear{2012}}

%
\begin{abstract}
In this paper, we show that reflecting Brownian motion in any bounded
domain $D$ can be approximated, as $k\to\infty$, by simple random
walks on ``maximal connected'' subsets of $(2^{-k}\Z^d) \cap D$ whose
filled-in interiors are inside of $D$.
\end{abstract}

%

\end{frontmatter}

\section{Introduction}\label{intro}

We proved in a recent article~\cite{BC} that reflecting Brownian
motion in a domain $D$ can be approximated by a sequence of random
walks on subsets $A_k$ of $(2^{-k}\Z^d) \cap D$. We chose $A_k$'s in a
``natural'' way, to be described in a moment.
Our main theorem in~\cite{BC} was limited to only some domains $D$
(``extension domains''). We also provided a counterexample showing that
random walks on $A_k$'s do not converge to the reflecting Brownian
motion in $D$ for some domains~$D$. In this paper, we will show in
Theorems~\ref{T34} and~\ref{T42}
that reflecting Brownian motion on any domain can be approximated
by a sequence of
discrete-time, as well as continuous-time,
random walks if the state spaces $D_k$ for the random walks are
constructed in a different ``natural'' way.

The sets $A_k$ were constructed in~\cite{BC} as follows.
First, we found the maximal connected set consisting of line segments
contained in $D$, joining neighboring vertices in $(2^{-k}\Z^d) \cap
D$. Then we let $A_k$ be the set of vertices in $(2^{-k}\Z^d) \cap D$
at the ends of these line segments. It turns out that the ``correct''
way (employed in the present article) to construct the state space for
the random walk is to start with the maximal connected set consisting
of cubes contained in $D$, with edge length $2^{-k}$ and vertices in
$(2^{-k}\Z^d) \cap D$. Then we let $D_k$ be the set of vertices in
$(2^{-k}\Z^d) \cap D$ which belongs to these cubes. Intuitively
speaking, $A_k$ may penetrate very thin crevices in $D$. The simple
random walk on $A_k$ may spend a nonnegligible amount of time in such
branches of $A_k$, but reflecting Brownian motion spends very little
time in sets with very small volume. Replacing edges in the
construction of $A_k$'s with cubes eliminates the mismatch between the
shapes of $D$ and the approximating discrete set.

The technical essence of the paper is Theorem~\ref{T21} which shows
that, in a sense, the Dirichlet
form for reflecting Brownian motion can be approximated
from below
by discrete Dirichlet forms. This theorem and the remaining part of the
proof of the main result are challenging because ``naive'' discrete
approximating schemes for
the Dirichlet form of reflecting Brownian motion
do not work; see Example~\ref{a161}.

In the rest of the \hyperref[intro]{Introduction}, we will review some basic facts about
reflecting Brownian motion in nonsmooth domains and elaborate on some
of the points mentioned above.

Reflecting Brownian motion in a bounded domain $D$ in $\R^d$
is a symmetric Markov process that behaves like Brownian motion
inside $D$ and is ``pushed'' back along the ``inward normal'' direction
at the boundary $\partial D$ of $D$.
It is a prototype of diffusions with boundary condition and can be used
to study heat equations with Neumann and Robin's boundary condition.
It is also widely used in modeling, for example,
in physics,
in queuing theory and in financial mathematics.
Reflecting Brownian motion has been studied by various authors
using various methods; see~\cite{BC,C}
and the references therein.
When $D$ is a bounded extension domain (see next paragraph
for its definition),
reflecting Brownian motion $X$ can be constructed as a strong Markov
process on $\overline D$
starting from every point in $\overline D$ except
a polar set.
Every bounded Lipschitz domain is an extension domain.
When $D$ is a general bounded domain, reflecting Brownian motion can
still be constructed
on $\overline D$, but typically it is no longer a strong Markov process.
In a recent paper~\cite{BC}, we developed three discrete approximation
schemes for reflecting Brownian motion in bounded domains, providing
effective ways to simulate the process in practice.
The first two approximation schemes are discrete-time
and continuous-time
simple random walks on grids $2^{-k} \Z^d\cap D$ inside $D$. For these
two approximation
schemes, we need to assume that $D$
is a bounded extension domain. A counter example is given in~\cite{BC},
showing that these approximation schemes do not work for some bounded domains.
However, the
third approximation scheme developed in~\cite{BC}, called myopic conditioning,
works for any bounded domain $D$.
Myopic conditioning generates a continuous-time and continuous-space
process and, therefore, it is not suited for computer simulations.
The purpose of this paper is to develop discrete-time
and continuous-time
simple random walk approximations on grids inside $D$ that
work for every
bounded domain $D$.

We now give a precise description of reflecting Brownian motion on
bounded domains.
Let $d\geq1$ and $D$ be a bounded domain in $\R^d$.
The Sobolev space
$W^{1,2}(D)$ of order $(1, 2)$ is\vadjust{\goodbreak} the space of $L^2(D)$-functions on
$D$ whose distributional derivative $\nabla f$ is also
$L^2(D)$-integrable. It is well known that $W^{1,2}(D)$ is a Hilbert
space under norm $\|f\|_{1,2}:= (\|f\|_{L^2(D)}^2+\|\nabla
f\|_{L^2(D)}^2 )^{1/2}$.
We define on $W^{1,2}(D)$ a bilinear form
\[
\EE(f, g)=\frac12 \int_D \nabla f(x) \cdot\nabla g(x) \,dx
\qquad\mbox{for } f, g\in W^{1,2}(D).
\]
It is known (see, e.g.,~\cite{CF}) that $(\EE, W^{1,2}(D))$ is a
Dirichlet form on $L^2(D; dx)$. When $C(\overline D)\cap W^{1,2}(D)$
is dense in both
$(C(\overline D), \| \cdot\|_\infty)$ and in
$(W^{1,2}(D),\break \| \cdot\|_{1,2})$, $(\EE, W^{1,2}(D))$ is a regular
Dirichlet form on $L^2(\overline D; m)$, where $m$ is the Lebesgue
measure on $D$ extended to $\overline D$ by setting $m(\partial
D)=0$. In this case, there is a continuous conservative strong
Markov process $X$ on $\overline D$ associated with $(\EE,
W^{1,2}(D))$, starting from quasi every point from $\overline D$.
The process is called (normally) reflecting Brownian motion on
$\overline D$.
It is known (see, e.g., Theorems 1 and 2 on pages 13 and~14 of
\cite{M}) that $(\EE, W^{1,2}(D))$ is a regular Dirichlet form
on $L^2(\overline D; m)$ if $D$ is
star-shaped with respect to a
point in $D$ or if $D$ has continuous boundary. Note that $(\EE,
W^{1,2}(\R^d))$ is a regular Dirichlet form on $L^2(\R^d; dx)$.
Hence $(\EE, W^{1,2}(D))$ is a regular Dirichlet form on
$L^2(\overline D; m)$ if $D$ is an
extension domain in the following sense: there is a linear
continuous operator
$T\dvtx W^{1,2}(D)\to W^{1,2}(\R^d)$
such that
$Tf=f$ a.e. on $D$ for every $f\in W^{1,2}(D)$.
Recall that a domain $D$ is called a \textit{locally uniform domain} if
there are
$\delta\in(0, \infty] $ and $C >0$ such that for every $x, y \in
D$ with $|x-y|<\delta$, there is a rectifiable curve $\gamma$ in $D$
connecting $x$ and
$y$ with $\operatorname{length}(\gamma) \leq C |x-y|$ and moreover,
\[
\min\bigl\{ |x-z|, |z-y| \bigr\} \leq C \operatorname{dist}\bigl(z, D^c\bigr)
\qquad\mbox{for every } z\in\gamma.
\]
A domain is said to be a \textit{uniform domain} if the above property
holds with
$\delta=\infty$.
The above definition is taken from V\"ais\"al\"a~\cite{Va}, where
various equivalent definitions are discussed.
Uniform domain and locally uniform domain are also called $(\eps,
\infty)$-domain
and $(\eps, \delta)$-domain, respectively, in~\cite{Jo}.
For example, the classical van Koch snowflake
domain in the conformal mapping theory is a uniform domain in
$\RR^2$. Note that every bounded Lipschitz domain is uniform, and
every \textit{nontangentially accessible domain} defined by Jerison
and Kenig in~\cite{JK} is a uniform domain (see (3.4) of
\cite{JK}),
while every Lipschitz domain is an $(\eps, \delta)$-domain.
It is proved in~\cite{Jo} that every locally uniform domain
is an extension domain.
However,
for general domain $D$, $(\EE, W^{1,2}(D), \EE)$
does not need
to be
regular on $L^2(\overline D; dx)$. A unit disk in $\R^2$ with a slit
removed is such an example. See page 14 of~\cite{M} for an example
of $D$ due to Kolsrud with $\partial D=\partial\overline D$ such
that the Dirichlet form $(\EE, W^{1,2}(D), \EE)$ is not regular on
$L^2(\overline D; dx)$. Nevertheless, for any domain $D\subset
\R^d$, one can always find a compact regularizing space $\widetilde D$ that
contains $D$ as a dense open subset such that $(\EE, W^{1,2}(D))$
becomes a regular Dirichlet space on $L^2(\widetilde D; \widetilde m)$,
where $\widetilde
m$ is the Lebesgue measure on $D$ extended to $\widetilde D$ by setting\vadjust{\goodbreak}
$\widetilde m (\widetilde D\setminus D)=0$; see~\cite{Fu} and~\cite{C}.
Let $\widetilde
X$ be the associated conservative strong Markov process on $\widetilde D$,
which can also be called reflecting Brownian motion on $D$. Let $X$
be the
projection
of $\widetilde X$ onto $\overline D$. Since for any given
time $t>0$, $\P_{\widetilde m}(\widetilde X_t\in\widetilde D \setminus
D)=0$, under the
normalized Lebesgue measure on $D$, $\widetilde X$ and $X$ have the same
finite-dimensional distributions.

A key technical element of this paper is to show that, for any bounded
domain $D$ in $\R^d$, there exists\vspace*{1pt} a sequence
$\{\varphi_j, j\geq1\}$ of bounded smooth functions on $D$ that is
dense in the Sobolev space $W^{1,2}(D)$, separates points in $D$ and
satisfies the property (\ref{e21}) described below. We can deduce from
its existence that there is a metric $\rho$ on $D$ (``refinement of the
Euclidean metric'') which induces the same Euclidean topology inside
$D$ and has the property that reflecting Brownian motion on $D$ can be
lifted as a strong Markov process on the $\rho$-closure $\widetilde D$
of $D$. This enables us to show that the random walk approximation on
grids whose filled-in interiors are inside $D$ works for reflecting
Brownian motions on arbitrary bounded domains. In this paper, we also
provide a proof that any weak limit of random walks on grids inside $D$
is a stationary symmetric Markov process (see Theorem~\ref{T32}). This
is a key step in proving that the weak limit is indeed the stationary
reflecting Brownian motion in $D$, using a Dirichlet form approach.
This claim was made in~\cite{BC} but regrettably no proof was given
there.

The rest of the paper is organized as follows.
In Section~\ref{S2}, we establish a result
(Theorem~\ref{T21})
regarding the Sobolev space $W^{1,2}(D)$ that will play
an important role in this paper. Though the result
is purely analytic, we employ some probabilistic
techniques in its proof. The proof that reflecting Brownian motion
in any bounded domain $D$ can be approximated by discrete-time
random walk on grids inside $D$ is given in Section~\ref{S3}.
The corresponding result for continuous-time random walk approximation
is presented in Section~\ref{S4}.

\section{Energy form estimates}\label{S2}

Let
$D\subset\R^d$ be a domain (connected open set) that has
finite Lebesgue measure. Fix an arbitrarily small $c_1\in(0,1)$
and a point $x_0\in D$. For each integer $k$, let $\A_k$
be the family of all closed $d$-dimensional cubes $Q\subset D$
with edge length $2^{-k}$, such that:
\begin{longlist}
\item the vertices of $Q$
belong to $(2^{-k}\Z)^d$;
\item the distance from $Q$ to $\prt
D$ is greater than $ c_1 2^{-k}$;
\item there exists a sequence of cubes $Q_1, Q_2, \ldots, Q_m=Q$,
satisfying (i) and~(ii), and such that $x_0\in Q_1$, and $Q_j \cap
Q_{j+1}$ is a $(d-1)$-dimensional cube, for all $j=1,2,\ldots,m-1$.
\end{longlist}

Since $D$ has a finite volume,
there is some $k_0\in\ZZ$ such that $\A_k = \varnothing$ for
every $k\leq k_0$.
Using scaling if necessary, we may and do assume
that $\A_k = \varnothing$ for $k\leq0$.
Let $D_k = \bigcup_{Q\in\A_k} Q$.
Let
$\A'_k$ be the family of all edges,
and let $\A''_k$ be the family of all vertices of all cubes $Q\in\A_k$.
We will write $\overline{xy}$ to denote the line segment with endpoints $x$
and $y$.
Note that if $\overline{xy} \in\A'_k$, then so is $\overline{yx}$.
Thus in the summation on the left-hand side of
(\ref{e21}), each line segment in $\A'_k$ is counted twice.\vadjust{\goodbreak}

\begin{theorem}\label{T21}
Suppose that $D\subset\R^d$ is domain with finite volume and
$c_1\in(0,1)$, $x_0\in D$ and $D_k$'s are defined as above. There
exists a countable sequence of bounded
functions $\{\vphi_j\}_{j\geq1} \subset W^{1,2}(D) \cap
C^\infty(D)$ such that:

\begin{longlist}
\item $\{\vphi_j\}_{j\geq1}$ is dense in
$W^{1,2}(D)$;

\item $\{\vphi_j\}_{j\geq1}$ separates points in
$D$;

\item for each $j\geq1$,
%
\begin{equation}
\label{e21} \limsup_{k\to\infty} 2^{k(2-d)} \sum
_{\overline{xy} \in\A'_k} \bigl(\vphi_j(x) - \vphi_j(y)
\bigr)^2 \leq2 \int_D \bigl|\nabla
\vphi_j (x)\bigr|^2 \,dx.
\end{equation}
\end{longlist}
\end{theorem}

\begin{pf}
\textit{Step} 1.
First note that the Sobolev space $(W^{1,2}(D), \| \cdot\|_{1,2})$
is separable. This can be seen as follows. Let $G_1$ be the
1-resolvent for the Dirichlet form $(\EE, W^{1,2}(D))$; that is,
$G_1$ is the linear operator from $L^2(D; m)$ to $W^{1,2}(D)$
uniquely defined by
\[
\EE_1 (G_1 f, g)= \int_D f(x)
g(x) \,dx \qquad\mbox{for every } g\in W^{1,2}(D).
\]
Here $\EE_1 (u, v):=\EE(u, v)+\int_D u(x) v(x) \,dx$. It follows
that $G_1 L^2(D; dx)$ is dense in the space $(W^{1,2}(D), \| \cdot\|_{1,2})$
and that
\[
\EE_1 (G_1f, G_1f)= \int
_D f G_1f(x) \,dx \leq\int_D
f(x)^2 \,dx.
\]
Since $L^2(D; dx)$ is separable, there is a sequence $\{f_k, k\geq
1\}$ of bounded functions that is dense in
$L^2(D; dx)$. Consequently,
$\{\eta^k:=G_1 f_k, k\geq1\}$ is a sequence of bounded functions
that is dense in
$(W^{1,2}(D), \| \cdot\|_{1,2})$.

Theorem 2 on page 251 of~\cite{E} implies that
for every function $ \eta^k$, there exists a sequence of functions
$\{ \eta^k_j, j\geq1\} \subset
W^{1,2}(D)\cap C^\infty(D)$ with the property that
$\lim_{j\to\infty} \|\eta^k_j - \eta^k \|_{1,2}=0$.
Moreover, the proof given in~\cite{E} shows that we can choose $
\eta^k_j$ so that
$\sup_{x\in D} | \eta^k_j (x)| \leq3 \sup_{x\in D} | \eta^k (x)|$.

\textit{Step} 2.
Constants $c_1, c_2, \ldots$ may change value from one ``step'' to
another in this proof.

We will use a regularized version of the distance function defined in
\cite{Stein}, Theorem 2, page 171. That theorem implies that there
exist $0<c_1,c_2,c_3, c_4<\infty$ such that for every integer $j$
there is a $C^\infty$ function $d_j\dvtx D \to(0,2^{-j}]$ with the
following properties:
%
\begin{eqnarray}
\label{a111}
c_1 \bigl(\dist(x, \prt D) \land2^{-j}\bigr)
&\leq& d_j(x) \leq c_2 \bigl(\dist(x, \prt D)
\land2^{-j}\bigr),
\\
\label{may141}
\sup_{x\in D} \bigl|\nabla d_j(x)\bigr| &\leq& c_3,
\\
\label{a112}
\sup_{x\in D} \biggl\llvert d_j (x) \,\frac\prt{\prt
x_i} \,\frac\prt{\prt x_m} d_j(x)\biggr
\rrvert&\leq& c_4 \qquad\mbox{for } 1\leq i,m\leq d.
\end{eqnarray}
By dividing $d_j$ by an appropriate constant, we may and will assume
from now on that (\ref{a111}) and (\ref{may141}) hold with $c_2 = c_3
= 1$.\vadjust{\goodbreak}

The existence of functions $d_j$ follows essentially from
\cite{Stein}, Theorem 2, pa\-ge~171.
The only difference between our claim and that in
\cite{Stein}, Theorem~2,
pa\-ge~171, is that~\cite{Stein} is concerned with the condition
\[
c_1 \dist(x, \prt D) \leq d_j(x) \leq c_2
\dist(x, \prt D),
\]
in place of our condition (\ref{a111}).
The method of proof given in~\cite{Stein} applies to (\ref{a111}) if
we subdivide cubes constructed in~\cite{Stein}, Section 1.2, page 167,
with edges longer than $2^{-j}$ into cubes with edge length $2^{-j}$.

\textit{Step} 3. Let $\psi\dvtx\R^d \to[0,\infty)$ be a $C^\infty$ ``mollifier''
with support in the ball $B(0,1/2)$ such that $\int_{B(0,1/2)} \psi
(x)\,dx =
1$. For $r>0$ let $\psi_r(x) = r^{-d}\psi( x/r)$, and note that $\sup_x
\psi_r(x) = c_1 r^{-d} $.

The function $\psi'_r(y):= \frac{\prt}{\prt r} \psi_{r}(y)$ is
$C^\infty$. It is supported in $B(0,r/2)$
and satisfies the condition $\int_{B(0,r/2)}
\psi'_{r}(y) \,dy = 0$.
Let $\| \cdot\|_1$ denote the $L^1$ norm with respect to the
Lebesgue measure restricted to $D$.
Note that $\|\psi'_r\|_1 = c_2 r^{-1}$ and $\|\psi'_r\lor0\|_1 = \|
\psi'_r\land0\|_1 = \|\psi'_r\|_1/2$.
Consider $x\in D$, and
let
\[
a^+_x( \cdot) = \frac{\psi'_{d_j(x)}( \cdot) \lor0}{\|\psi
'_{d_j(x)}\lor0\|_1} \quad\mbox{and}\quad a^-_x(
\cdot) = - \frac
{\psi'_{d_j(x)}( \cdot) \land0} {
\|\psi'_{d_j(x)}\land0\|_1}.
\]
The functions\vspace*{1pt} $a^+_x( \cdot)$ and $a^-_x( \cdot)$ are
probability density functions. Let $A^+_x $ and $A^-_x$ be independent
$\RR^d$-valued random variables
with densities $a^+_x( \cdot)$ and $a^-_x( \cdot)$,
respectively. Let
$\int_{A^+_x}^{A^-_x} d \xi$
denote the integral with respect to the length measure on the line
segment joining $A^+_x $ and $A^-_x$.
Clearly, the measure
$\E\int_{A^+_x}^{A^-_x} d\xi$
is supported on
$\overline{B(0, d_j(x)/2)}$. We
will now show that it
has a density bounded above by $c_3 d_j(x)^{1-d}$.
In other words, for every set $K\subset D$,
\[
\E\int_{A^+_x}^{A^-_x} \bone_K(\xi) \,d \xi
\leq c_3 d_j(x)^{1-d} m \bigl(K\cap B\bigl(0,
d_j(x)/2\bigr)\bigr).
\]

Clearly the functions
$a^+_x( \cdot)$ and $a^-_x( \cdot)$ are bounded by $\alpha_x
:= c_4 d_j(x)^{-d}$. Consider any $z\in B(x,d_j(x)/2)$ and small
$\delta>0$.
The probability that at least one of the random points $A^+_x $ or
$A^-_x$ belongs to $B(z, 2\delta)$ is less than $c_5 \delta^d\alpha
_x\leq c_6 \delta^d d_j(x)^{-d}$.

For $k\geq2$, the probability of the event $F_k = \{A^+_x \in
B(z,\delta2^k) \setminus B(z,\delta2^{k-1})\}$ is bounded by $c_7
\delta^d 2^{dk} \alpha_x$.
Let $G$ be the intersection of $B(x,d_j(x)/2)$ and the smallest cone
with vertex $A^+_x$ containing $B(z, \delta)$. The conditional
probability, given~$F_k$, of the event $\{A^-_x\in G\}$
is bounded by $\alpha_x$ times the volume of $G$; hence, it is bounded
by $c_8 \alpha_x d_j(x) (2^{-k} d_j(x))^{d-1} = c_9 \alpha_x d_j(x)^d
2^{k(1-d)}$.
Multiplying the two estimates and
summing over $k\geq2$, such that $B(z,\delta2^{k-1}) $ does not
contain $B(x,d_j(x)/2)$, gives
the bound
\[
\sum_{\delta2^k \leq2 d_j(x)} c_7 \delta^d
2^{dk} \alpha_x c_9 \alpha_x
d_j(x)^d 2^{k(1-d)} \leq c_{10}
\delta^{d-1} d_j(x)^{1-d}.
\]
Adding to this quantity $c_6 \delta^d d_j(x)^{-d}$ [the estimate
representing the case when
at least one of the random points $A^+_x $ or $A^-_x$ belongs to $B(z,
2\delta)$] gives a similar bound
$c_{11} \delta^{d-1} d_j(x)^{1-d}$. The last quantity is an upper
bound for the probability that the line segment joining $A^+_x $ and
$A^-_x$ intersects $B(z, 2\delta)$.
Since $\int_{A^+_x}^{A^-_x} \bone_{B(z, 2\delta)}(\xi) \,d\xi\leq2
\delta$ with probability 1, we obtain
\[
\E\int_{A^+_x}^{A^-_x} \bone_{B(z, 2\delta)}(\xi) \,d\xi
\leq c_{11} \delta^{d-1} d_j(x)^{1-d}2
\delta= c_{12} \delta^{d} d_j(x)^{1-d}.
\]
This estimate holds for all $z\in B(x, d_j(x)/2)$ and all small $\delta
>0$ so the density of the measure
$\E\int_{A^+_x}^{A^-_x} d\xi$
is bounded by $c_{13} d_j(x)^{1-d} $.

We will also need the following version of the above estimate.
Let $\psi''_r(y) = \frac{\prt^2}{\prt r^2} \psi_r(y)$, $b^+_x(
\cdot) = (\psi''_{d_j(x)}( \cdot) \lor0)/\|\psi''_{d_j(x)}\lor0\|_1$ and
$b^-_x( \cdot) = -(\psi''_{d_j(x)}( \cdot) \land0)/\break\|\psi
''_{d_j(x)}\land0\|_1$.
Note that $\|\psi'_r\|_1 = c_{14} r^{-2}$,
$\int_{B(0, r/2)} \psi''_r (y) \,dy=0$ and so
$\|\psi''_r\lor0\|_1 = \|\psi''_r\land0\|_1 = \|\psi''_r\|_1/2$.
The functions $b^+_x( \cdot)$ and $b^-_x( \cdot)$ are
probability density functions. Let $B^+_x $ and $B^-_x$
be independent $\RR^d$-valued
random variables with densities $b^+_x( \cdot)$ and $b^-_x( \cdot
)$.
The measure $\E\int_{B^+_x}^{B^-_x} d\xi$ has a density bounded
above by $c_{15} d_j(x)^{1-d}$.
In other words, for every set $K\subset D$,
\[
\E\int_{B^+_x}^{B^-_x} \bone_K(\xi) \,d\xi
\leq c_{15} d_j(x)^{1-d} m \bigl(K\cap B\bigl(0,
d_j(x)/2\bigr)\bigr).
\]
We omit the proof because it is analogous to the one given above.

\textit{Step} 4. Consider a function $ \eta\in W^{1,2}(D) \cap C^\infty
(D)$ and for integer
$j\geq1$ and $x\in D$, let
%
\begin{equation}
\label{a91} \eta_j (x) = \int_{B(0, d_j(x)/2)}
\psi_{d_j(x)}(y) \eta(x-y) \,dy.
\end{equation}

We will show that $\eta_j \in W^{1,2}(D)\cap C^\infty(D)$
and $\eta_j \to\eta$ in
$W^{1,2}(D)$ as $j\to\infty$.
Since $\eta$ and $d_j$ are $C^\infty$ functions,
so is $\eta_j$.

We have
%
\begin{eqnarray}
\label{a92} \eta_j(x)^2 &=& \biggl(\int
_{B(0, d_j(x)/2)} \psi_{d_j(x)}(y) \eta(x-y) \,dy
\biggr)^2
\nonumber\\
&\leq&\int_{B(0, d_j(x)/2)} \psi_{d_j(x)}(y)^2 \,dy \int
_{B(0, d_j(x)/2)} \eta(x-y)^2 \,dy
\nonumber\\[-8pt]\\[-8pt]
&\leq& c_1 \bigl(d_j(x)^{-d}
\bigr)^2 d_j(x)^d \int_{B(0, d_j(x)/2)}
\eta(x-y)^2 \,dy
\nonumber
\\
&\leq& c_2 d_j(x)^{-d} \int
_{B(0, d_j(x)/2)} \eta(x-y)^2 \,dy.
\nonumber
\end{eqnarray}
Suppose that $z\in B(x, d_j(x)/2)$.
Then
\[
\dist(z,\prt D) \geq\dist(x,\prt D) - |x-z| \geq d_j(x) - |x-z|
\geq d_j(x)/2.
\]
Hence
\[
d_j(z) \geq c_3\bigl(\dist(z, \prt D)
\land2^{-j}\bigr) \geq c_3\bigl(d_j(x)/2
\land2^{-j}\bigr) = c_3 d_j(x)/2.
\]
Therefore, $c_4 d_j(z) \geq d_j(x)/2$ and $\bone_{B(x, d_j(x)/2)}(z)
\leq\bone_{B(z,c_4 d_j(z))}(x)$.
Assuming again that $z\in B(x, d_j(x)/2)$, we obtain
\begin{eqnarray*}
\dist(z,\prt D) &\leq&\dist(x,\prt D) + |x-z| \leq\dist(x,\prt D) +
d_j(x)/2
\\[-2pt]
&\leq&\dist(x,\prt D) + \dist(x,\prt D)/2 < 2\dist(x,\prt D).
\end{eqnarray*}
Hence
\begin{eqnarray*}
d_j(x) &\geq& c_5\bigl( \dist(x,\prt D)
\land2^{-j}\bigr) \geq c_5\bigl( \dist(z,\prt D)/2
\land2^{-j}\bigr) \\[-2pt]
&\geq& c_5\bigl(d_j(z)/2
\land2^{-j}\bigr) = c_5 d_j(z)/2.
\end{eqnarray*}
This implies that
%
\begin{equation}
\label{a171} d_j(x)^{-d} \bone_{B(x, d_j(x)/2)}(z) \leq
c_5^{-d}\bigl(d_j(z) /2\bigr)^{-d}
\bone_{B(z,c_4 d_j(z))}(x).
\end{equation}
For later reference we derive an inequality that is slightly more
general than what is needed in this step. For a set $Q\subset D$, let
$\widehat Q = \bigcup_{x\in Q} B(x, d_j(x)/2)$.
We combine (\ref{a92}) and (\ref{a171}) to see that
%
\begin{eqnarray}
\label{a1420} \int_Q \eta_j(x)^2
\,dx &\leq& c_2 \int_Q d_j(x)^{-d}
\int_{B(0, d_j(x)/2)} \eta(x-y)^2 \,dy \,dx
\nonumber\\[-2pt]
&=& c_2 \int_Q d_j(x)^{-d}
\int_{B(x, d_j(x)/2)} \eta(z)^2 \,dz \,dx
\nonumber
\\[-2pt]
&=& c_2 \int_Q \int_{\widehat Q}
d_j(x)^{-d} \bone_{B(x, d_j(x)/2)}(z) \eta(z)^2
\,dz \,dx
\\[-2pt]
&\leq& c_2 \int_{\widehat Q} \int_Q
c_5^{-d}\bigl(d_j(z) /2\bigr)^{-d}
\bone_{B(z, c_4 d_j(z))}(x)\,dx\, \eta(z)^2 \,dz
\nonumber
\\[-2pt]
& \leq& c_6 \int_{\widehat Q} \eta(z)^2 \,dz.
\nonumber
\end{eqnarray}
In particular, the inequality applies to $Q = D = \widehat Q$. Hence
%
\begin{equation}
\label{a102} \int_D \eta_j(x)^2
\,dx \leq c_6 \int_D \eta(z)^2 \,dz.
\end{equation}

For any $x\in D$, $j$ and $1\leq i \leq d$,
note that $\psi_{d_j(x)}(y)=0$ for $y\notin B(0, d_j(x)/2)$,
and $\eta(x-y)$ is differentiable in
$y\in\overline{B(0, d_j(x)/2)}$.
So we have
%
\begin{eqnarray}
\label{a101}
\biggl(\frac{\prt}{\prt x_i} \eta_j(x)
\biggr)^2 &=& \biggl(\frac{\prt}{\prt x_i} \int_{B(0, d_j(x)/2)}
\psi_{d_j(x)}(y) \eta(x-y) \,dy \biggr)^2
\nonumber\\[-2pt]
&=& \biggl( \int_{B(0, d_j(x)/2)} \psi_{d_j(x)}(y)\,
\frac{\prt}{\prt
x_i} \eta(x-y) \,dy\nonumber\\[-2pt]
&&\hspace*{21pt}{} +\int_{\R^d} \biggl(
\frac{\prt}{\prt x_i} \psi_{d_j(x)}(y) \biggr) \eta(x-y) \,dy
\biggr)^2
\nonumber
\\[-2pt]
&=& \biggl( \int_{B(0, d_j(x)/2)} \psi_{d_j(x)}(y)
\,\frac{\prt}{\prt
x_i} \eta(x-y) \,dy \nonumber\\[-2pt]
&&\hspace*{5pt}{} +\int_{B(0, d_j(x)/2)} \biggl(
\frac{\prt}{\prt x_i} \psi_{d_j(x)}(y) \biggr) \eta(x-y) \,dy
\biggr)^2
\\[-2pt]
&\leq& 2 \biggl( \int_{B(0, d_j(x)/2)} \psi_{d_j(x)}(y)
\,\frac{\prt}{\prt
x_i} \eta(x-y) \,dy \biggr)^2\nonumber\\[-2pt]
&&{} +2 \biggl( \int
_{B(0, d_j(x)/2)} \biggl(\frac{\prt}{\prt x_i} \psi_{d_j(x)}(y) \biggr)
\eta(x-y) \,dy \biggr)^2
\nonumber
\\[-2pt]
&\leq& 2\int_{B(0, d_j(x)/2)} \psi_{d_j(x)}(y)^2 \,dy
\int_{B(0, d_j(x)/2)} \biggl( \frac{\prt}{\prt x_i}\eta(x-y)
\biggr)^2 \,dy
\nonumber
\\[-2pt]
&&{} + 2 \biggl( \int_{B(0, d_j(x)/2)} \biggl(\frac{\prt}{\prt x_i}
\psi_{d_j(x)}(y) \biggr) \eta(x-y) \,dy \biggr)^2
\nonumber
\\[-2pt]
&\leq& c_7 d_j(x)^{-d} \int
_{B(0, d_j(x)/2)} \biggl( \frac{\prt}{\prt x_i}\eta(x-y)
\biggr)^2 \,dy\nonumber\\[-2pt]
&&{} + 2 \biggl( \int_{B(0, d_j(x)/2)} \biggl(
\frac{\prt}{\prt x_i} \psi_{d_j(x)}(y) \biggr) \eta(x-y) \,dy
\biggr)^2.
\nonumber
\end{eqnarray}

Recall that $\psi'_r(y):= \frac{\prt}{\prt r} \psi_{r}(y)$ is a
$C^\infty$ function supported in $B(0,1/2)$ with $\int_{B(0,1/2)}
\psi'_{r}(y) \,dy = 0$.
It follows from
(\ref{may141}) that
${\sum_{i=1}^d} |{\frac{\prt}{\prt x_i} d_j(x)}| \leq d$.
We have
%
\begin{eqnarray}
\label{a141} \sum_{i=1}^d \biggl\llvert
\frac{\prt}{\prt x_i} \psi_{d_j(x)}(y)\biggr\rrvert&=& \sum
_{i=1}^d \biggl\llvert\frac{\prt}{\prt x_i}
d_j(x)\biggr\rrvert\bigl|\psi'_{d_j(x)}(y)\bigr| \nonumber\\[-9pt]\\[-9pt]
&\leq& d \bigl|
\psi'_{d_j(x)}(y)\bigr| \leq c_9 d_j(x)^{-d-1}.\nonumber
\end{eqnarray}
Recall\vspace*{1pt} the definitions of the random variables $A_x^+$ and $A_x^-$
from step 3. It follows from
$\int_{B(0, d_j(x)/2)} \psi'_{d_j(x)}(y) \,dy=0$
and $\int_{B(0, d_j(x)/2)} |\psi'_{d_j(x)}(y)| \,dy\leq c_{10}
/d_j(x)$ that
%
\begin{eqnarray}
\label{a146}
&&
\biggl\llvert\int_{B(0, d_j(x)/2)} \biggl(
\frac{\prt}{\prt x_i} \psi_{d_j(x)}(y) \biggr) \eta(x-y) \,dy \biggr
\rrvert\nonumber\\[-2pt]
&&\qquad\leq
c_{10} d_j(x)^{-1} \bigl\llvert\E\bigl(\eta
\bigl(x-A^+_x\bigr) - \eta\bigl(x-A^-_x\bigr)\bigr)
\bigr\rrvert
\nonumber\\[-8pt]\\[-8pt]
&&\qquad\leq c_{10} d_j(x)^{-1}\E\bigl\llvert\eta
\bigl(x-A^+_x\bigr) - \eta\bigl(x-A^-_x\bigr) \bigr
\rrvert\nonumber\\
&&\qquad\leq c_{10} d_j(x)^{-1}\E\int
_{A^+_x}^{A^-_x} \bigl|\nabla\eta(x-z)\bigr| \,dz,
\nonumber
\end{eqnarray}
where the last integral is along a line segment from
$A^+_x$ to $A^-_x$.
By step 3, the measure $\E\int_{A^+_x}^{A^-_x} dz$ has a density
that is bounded above by $c_{11} d_j(x)^{1-d}$
and vanishes outside of the ball $B(0, d_j(x)/2)$.
In other words, for every set $K\subset D$, $\E\int_{A^+_x}^{A^-_x}
\bone_K(z) \,dz\leq c_{11} d_j(x)^{1-d}
m(K\cap B(0, d_j(x)/2))$.
It follows that
%
\begin{eqnarray}
\label{a143}
&&
\biggl\llvert\int_{B(0, d_j(x)/2)} \biggl(
\frac{\prt}{\prt x_i} \psi_{d_j(x)}(y) \biggr) \eta(x-y) \,dy \biggr
\rrvert\nonumber\\
&&\qquad\leq
c_{10} d_j(x)^{-1}\E\int_{A^+_x}^{A^-_x}
\bigl|\nabla\eta(x-z)\bigr| \,dz
\nonumber\\[-8pt]\\[-8pt]
&&\qquad\leq c_{10} d_j(x)^{-1} \int
_{B(0, d_j(x)/2)} \bigl|\nabla\eta(x-z)\bigr| c_{11}
d_j(x)^{1-d} \,dz \nonumber\\
&&\qquad= c_{12} d_j(x)^{-d}
\int_{B(0, d_j(x)/2)} \bigl|\nabla\eta(x-z)\bigr| \,dz.
\nonumber
\end{eqnarray}
This implies that
\begin{eqnarray*}
&&
\biggl( \int_{B(0, d_j(x)/2)} \biggl(\frac{\prt}{\prt x_i}
\psi_{d_j(x)}(y) \biggr) \eta(x-y) \,dy \biggr)^2 \\
&&\qquad\leq \biggl(
c_{12} d_j(x)^{-d} \int_{B(0, d_j(x)/2)}
\bigl|\nabla\eta(x-z)\bigr| \,dz \biggr)^2
\\
&&\qquad\leq c_{13} d_j(x)^{-2d} d_j(x)^{d}
\int_{B(0, d_j(x)/2)} \bigl|\nabla\eta(x-z)\bigr|^2 \,dz \\
&&\qquad=
c_{13} d_j(x)^{-d} \int_{B(0, d_j(x)/2)}
\bigl|\nabla\eta(x-z)\bigr|^2 \,dz.
\end{eqnarray*}
We combine this estimate with (\ref{a101}) to obtain
\begin{eqnarray*}
&& \biggl(\frac{\prt}{\prt x_i} \eta_j(x) \biggr)^2
\\
&&\qquad\leq c_7 d_j(x)^{-d} \int
_{B(0, d_j(x)/2)} \biggl( \frac{\prt}{\prt x_i}\eta(x-y)
\biggr)^2 \,dy \\
&&\qquad\quad{}+ c_{14} d_j(x)^{-d}
\int_{B(0, d_j(x)/2)} \bigl|\nabla\eta(x-z)\bigr|^2 \,dz.
\end{eqnarray*}
Summing over $i$ yields
\[
\bigl\llvert\nabla\eta_j(x)\bigr\rrvert^2 \leq
c_{15} d_j(x)^{-d} \int_{B(0, d_j(x)/2)}
\bigl|\nabla\eta(x-z)\bigr|^2 \,dz.
\]
Recall that we write $\widehat Q = \bigcup_{x\in Q} B(x, d_j(x)/2)$ for
$Q\subset D$.
The same argument that leads from (\ref{a92}) to (\ref{a102}) gives
%
\begin{equation}
\label{e214} \int_Q \bigl\llvert\nabla
\eta_j(x)\bigr\rrvert^2 \,dx \leq c_{16} \int
_{\widehat Q} \bigl|\nabla\eta(x)\bigr|^2 \,dx.
\end{equation}
This formula and (\ref{a102}) show that $\eta_j \in W^{1,2}(D)$. We
have pointed out earlier in the proof that $\eta_j \in C^\infty(D)$.

Let $K_\eps= \{x\in D\dvtx\dist(x, D^c) >\eps\}$.
We will show that $\eta_j \to\eta$ in
$W^{1,2}(D)$ as $j\to\infty$.
To see this, fix an arbitrarily
small $\delta>0$ and find $\eps>0$ so small that
%
\begin{equation}
\label{e215} \int_{K_{2\eps}^c} \bigl(\eta(x)^2 + \bigl|
\nabla\eta(x)\bigr|^2 \bigr) \,dx < \delta.
\end{equation}
Note that the integral in the above formula is over the set ${K_{2\eps
}^c}$, not ${K_{\eps}^c}$.
Since $\overline{K}_{\eps}
\subset D$ and $\eta$ is $C^\infty$, we have
%
\begin{equation}
\label{e215b} \lim_{j\to\infty} \int_{K_\eps} \bigl(\bigl(
\eta_j(x)-\eta(x)\bigr)^2 + \bigl(\bigl|\nabla
\eta_j(x)\bigr|^2 - \bigl|\nabla\eta(x)\bigr|^2\bigr) \bigr)
\,dx = 0,
\end{equation}
because the integrand converges to 0 uniformly. It suffices to show
that there exists a constant $c_{17} < \infty$, not depending on
$\delta$ or $\eps$, such that for large
$j$, we have
%
\begin{equation}
\label{a104} \int_{K_\eps^c} \bigl(\eta_j(x)^2
+ \bigl|\nabla\eta_j(x)\bigr|^2 \bigr) \,dx < c_{17}
\delta.
\end{equation}
By (\ref{a1420}) applied to $Q = K_\eps^c$,
\[
\int_{K_\eps^c} \eta_j(x)^2 \,dx \leq
c_6 \int_{K_{2\eps}^c} \eta(z)^2 \,dz \leq
c_6 \delta,
\]
while (\ref{e214}) and (\ref{e215}) imply that
\[
\int_{K_{\eps}^c} \bigl|\nabla\eta_j(x)\bigr|^2 \,dx
\leq c_{16} \delta.
\]
The last two estimates yield (\ref{a104}) and complete the proof of
the claim that
$\eta_j \to\eta$ in
$W^{1,2}(D)$ as $j\to\infty$.

\textit{Step} 5.
Recall the constant $c_1\in(0,1)$ and sets
$\{D_k, k\geq1\}$
from the beginning of this section. For each integer $k\geq1$, let $\B_k$
be the family of all closed $d$-dimensional cubes $Q\subset D$
with edge length $2^{-k}$, such that: (i) the vertices of $Q$
belong to $(2^{-k}\Z)^d$; (ii) the distance from $Q$ to $\prt
D$ is greater than $ c_1 2^{-k}$; and (iii) the interiors of $Q$ and
$D_k$
are disjoint.
Let $\M_1 = \B_1$ and let $\M_k\subset\B_k$
consist of those cubes in $\B_k$ that are not a subset of any cube
in $ \B_{k-1}$ for $k\geq2$.\vadjust{\goodbreak}

We recall that for a set $Q\subset D$, $\widehat Q = \bigcup_{x\in Q} B(x,
d_j(x)/2)$.
We claim that there exists $M<\infty$,
independent of $j$, such that every point $x\in D$
belongs to at most $M$ distinct sets of the form $\widehat Q$ where
$Q\in
\bigcup_k \M_k$.
This claim can be proved in a way that is totally analogous
to the proof of~\cite{Stein}, Proposition~3, page 169, so we omit
its proof.

\textit{Step} 6. We have shown in step 1 that we can find a
sequence of bounded functions $\{\eta^k\}_{k\geq1} $ in $
W^{1,2}(D) \cap C^\infty(D)$ that is dense in $W^{1,2}(D)$.
Let $\{\eta^k_j\}_{j\geq1}$ be a sequence constructed from
$\eta^k$ as in (\ref{a91}). The family $\{\eta^k_j\}_{k,j\geq1}$
is dense in $W^{1,2}(D)$ and consists of bounded $C^\infty$
functions. Let us relabel the family $\{\eta^k_j\}_{k,j\geq1}$
as $\{\vphi_j\}_{j\geq1}$. We see that the family
$\{\vphi_j\}_{j\geq1}$ consists of bounded functions in
$W^{1,2}(D) \cap C^\infty(D)$ and satisfies part (i) of the
theorem.

By adding an appropriate sequence of functions in
$C^\infty_b(\overline D)$, if necessary, we can assume that
condition (ii) is satisfied by $\{\vphi_j\}_{ j\geq1}$.

We will show that (\ref{e21}) holds for
$\vphi_j$ for each fixed $j\geq1$.
Some functions $\vphi_j$ belong to $C^\infty_b(\overline D)$. It is
easy to see that
(\ref{e21}) holds for such functions. Hence, we will assume that
$\vphi_j$ belongs to the
family $\{\eta^k_j\}_{k,j\geq1} $. Then there exists a
function $\vphi$ in $W^{1,2}(D) \cap C^\infty(D)$ such that $\vphi_j$
was constructed from $\vphi$ as in (\ref{a91}).

Fix an arbitrarily small
$\eps>0$ and find an integer $R$ so large that
\[
\| \nabla\vphi
\bone_{D_R^c} \|_{L^2(D)} < \eps.
\]
Note that $D_R$ is a compact set
and
choose an integer $S>R$ so large that $D_R$ is a subset of the
interior of $D_S$.
Recall
$\A_k,\A'_k$ and $\A''_k$
from the beginning of this section.
Let $e_i$ denote the unit vector in the positive direction of $x_i$-axis.
Since $\vphi_j$ is
$C^\infty(D)$, we have by the
mean value theorem for some
$\theta^+_i(x) \in\overline{x, x + 2^{-k}e_i }$
and $\theta^-_i(x) \in\overline{x, x - 2^{-k}e_i }$ that
\begin{eqnarray*}
&&
\limsup_{k\to\infty} 2^{k(2-d)} \sum_{\overline{x,y} \in\A'_k,
\overline{xy} \subset D_S }
\bigl(\vphi_j(x) - \vphi_j(y)\bigr)^2
\\
&&\qquad = \limsup_{k\to\infty} 2^{k(2-d)} \sum_{x\in\A''_k\cap D_S}
\sum_{y\in\A''_k\cap D_S\dvtx \overline{xy}\in A'_k} \bigl(\vphi_j(x) -
\vphi_j(y)\bigr)^2
\\
&&\qquad\leq\limsup_{k\to\infty} 2^{k(2-d)} \sum_{x\in\A''_k\cap D_S}
\sum_{i=1}^d \bigl( \bigl(
\vphi_j(x) - \vphi_j\bigl(x + 2^{-k}e_i
\bigr)\bigr)^2\\
&&\qquad\quad\hspace*{123pt}{} + \bigl(\vphi_j(x) - \vphi_j
\bigl(x - 2^{-k}e_i\bigr)\bigr)^2 \bigr)
\\
&&\qquad= \limsup_{k\to\infty} 2^{k(2-d)} \sum_{x \in\A''_k\cap D_S}
\sum_{i=1}^d \biggl( \biggl\llvert
\frac{\partial\vphi_j (\theta^+_i (x ))} {\partial x_i} \biggr\rrvert^2
+ \biggl\llvert\frac{\partial\vphi_j (\theta^+_i (x ))} {
\partial x_i}
\biggr\rrvert^2 \biggr) 2^{-2k}
\cr
&&\qquad\leq2 \int
_{D} |\nabla\vphi_j|^2.
\end{eqnarray*}

It suffices to find $c_1<\infty$ independent of $\eps, S$ and $R$ and
such that
\[
\limsup_{k\to\infty} 2^{k(2-d)} \sum_{\overline{xy} \in\A'_k, \overline
{xy} \not\subset D_S}
\bigl(\vphi_j(x) - \vphi_j(y)\bigr)^2 \leq
c_1 \eps.
\]
We assumed that $D_R$ is a subset of the
interior of $D_S$, so for large $k$, if
$\overline{x,y} \in\A'_k$ and $ \overline{xy} \not\subset D_S$,
then $ \overline{xy} \subset D_R^c$. Hence, it will suffice to
find $c_1< \infty$ such that
%
\begin{equation}
\label{eA83} \limsup_{k\to\infty} 2^{k(2-d)} \sum
_{\overline{xy} \in\A'_k, \overline{xy} \subset D_R^c} \bigl(\vphi_j(x)
- \vphi_j(y)
\bigr)^2 \leq c_1 \eps.
\end{equation}

Recall the notation from step 5.
Consider a large integer $k$, $\ell\leq k$ and $Q\in\M_\ell$.
Suppose that
%
\begin{equation}
\label{e216} \sum_{x,y \in(2^{-k}\Z)^d \cap Q, |x- y|=2^{-k}} \bigl(
\vphi_j(x) - \vphi_j(y)\bigr)^2 = a.
\end{equation}
Let $N = 2^{(k-\ell)d}$, and let $\{Q_1, Q_2, \ldots, Q_N\}$
be the family of all cubes such that
$Q_n \in\B_k$ and $Q_n \subset Q$.
Let $a_n $ be the maximum\vspace*{1pt} of $(\vphi_j(x) - \vphi_j(y))^2$ over all
pairs $x,y \in(2^{-k}\Z)^d \cap Q_n$ such that $|x- y|=2^{-k}$.
By the mean value theorem, there is some $z$ in the line segment
joining $x$ and $y$ in $Q_n$ such that
%
\begin{equation}
\label{eA43} \bigl|\nabla\vphi_j (z)\bigr| \geq a_n^{1/2}
2^{k}.
\end{equation}

It is easy to check that
%
\begin{equation}
\label{a125} d2^{d-1}\sum_{n=1}^N
a_n \geq a.
\end{equation}

\textit{Step} 7.
In this step, we will prove (\ref{eA83}).
We have
%
\begin{eqnarray}
\label{a172}
\frac{\prt}{\prt x_i} \biggl(\frac{\prt}{\prt x_m} \vphi_j(x)
\biggr)
&=&\frac{\prt}{\prt x_i} \biggl(\frac{\prt}{\prt x_m} \biggl( \int
_{\R^d} \psi_{d_j(x)}(y) \vphi(x-y) \,dy \biggr) \biggr)
\nonumber\\
&=& \frac{\prt}{\prt x_i} \biggl( \int_{\R^d}
\psi_{d_j(x)}(y) \,\frac{\prt}{\prt x_m} \vphi(x-y) \,dy\nonumber\\
&&\hspace*{22.3pt}{} +\int_{\R^d}
\biggl(\frac{\prt}{\prt x_m} \psi_{d_j(x)}(y) \biggr) \vphi(x-y) \,dy
\biggr)
\nonumber
\\
&=& \frac{\prt}{\prt x_i} \biggl( \int_{\R^d}
\psi_{d_j(x)}(x-y) \,\frac{\prt}{\prt x_m} \vphi(y) \,dy\nonumber\\
&&\hspace*{23pt}{} +\int_{\R^d}
\biggl(\frac{\prt}{\prt x_m} \psi_{d_j(x)}(y) \biggr) \vphi(x-y) \,dy
\biggr)
\nonumber\\[-8pt]\\[-8pt]
&=& \int_{\R^d} \biggl( \frac{\prt}{\prt x_i}
\psi_{d_j(x)}(x-y) \biggr) \,\frac{\prt}{\prt x_m} \vphi(y) \,dy
\nonumber\\
&&{} +\int_{\R^d} \biggl(\frac{\prt}{\prt x_i}\,\frac{\prt}{\prt x_m}
\psi_{d_j(x)}(y) \biggr) \vphi(x-y) \,dy\nonumber\\
&&{} +\int_{\R^d}
\biggl(\frac{\prt}{\prt x_m} \psi_{d_j(x)}(y) \biggr) \,\frac{\prt}{\prt x_i}
\vphi(x-y) \,dy
\nonumber
\\
&=& 2\int_{B(x, d_j(x)/2)} \biggl( \frac{\prt}{\prt x_i}
\psi_{d_j(x)}(x-y) \biggr) \,\frac{\prt}{\prt x_m} \vphi(y) \,dy
\nonumber
\\
&&{} +\int_{B(0, d_j(x)/2)} \biggl(\frac{\prt}{\prt x_i}\,\frac
{\prt}{\prt x_m}
\psi_{d_j(x)}(y) \biggr) \vphi(x-y) \,dy.
\nonumber
\end{eqnarray}
We estimate the first of the last two integrals, using (\ref{a141})
as follows:
%
\begin{eqnarray}
\label{a142}
&&\biggl\llvert\int_{B(x, d_j(x)/2)} \biggl(
\frac{\prt}{\prt x_i} \psi_{d_j(x)}(x-y) \biggr) \,\frac{\prt}{\prt x_m}
\vphi(y) \,dy
\biggr\rrvert
\nonumber
\\
&&\qquad\leq \int_{B(x, d_j(x)/2)} c_1 d_j(x)^{-d-1}
\biggl\llvert\frac{\prt}{\prt
x_m} \vphi(y)\biggr\rrvert \,dy \\
&&\qquad\leq c_2
d_j(x)^{-d-1} \int_{B(x, d_j(x)/2)} \bigl|\nabla
\vphi(z)\bigr| \,dz.\nonumber
\end{eqnarray}
To estimate the second integral, we apply the same method as in the
derivation of (\ref{a143}).
Recall that $\psi'_r(y) = \frac{\prt}{\prt r} \psi_{r}(y)$ is a
$C^\infty$ function supported in $B(0,1/2)$ with $\int_{B(0,1/2)}
\psi'_{r}(y) \,dy = 0$.
The function $\psi''_r(y) = \frac{\prt^2}{\prt r^2} \psi_{r}(y)$
is $C^\infty$. It is supported in $B(0,1/2)$ and satisfies the
condition $\int_{B(0,1/2)} \psi''_{r}(y) \,dy = 0$.
Note that $\|\psi'_r\|_1 = c_3 r^{-1}$, $\|\psi'_r\lor0\|_1 = \|\psi
'_r\land0\|_1 = \|\psi'_r\|_1/2$, $\|\psi''_r\|_1 = c_4 r^{-2}$, $\|
\psi''_r\lor0\|_1 = \|\psi''_r\land0\|_1 = \|\psi''_r\|_1/2$.
We have
\begin{eqnarray*}
\frac{\prt}{\prt x_i} \,\frac{\prt}{\prt x_m} \psi_{d_j(x)}(y) &=& \frac
{\prt}{\prt x_i}
\biggl( \biggl( \frac{\prt}{\prt x_m} d_j(x) \biggr)
\psi'_{d_j(x)}(y) \biggr)
\\
&=& \biggl(\frac{\prt}{\prt x_i} \,\frac{\prt}{\prt x_m} d_j(x) \biggr)
\psi'_{d_j(x)}(y)\\
&&{} + \biggl( \frac{\prt}{\prt x_m}
d_j(x) \biggr) \biggl( \frac{\prt}{\prt x_i} d_j(x) \biggr)
\psi''_{d_j(x)}(y).
\nonumber
\end{eqnarray*}
This implies that
%
\begin{eqnarray}
\label{a151}\qquad
&&
\biggl\llvert\int_{B(0, d_j(x)/2)} \biggl(
\frac{\prt}{\prt x_i}\,\frac
{\prt}{\prt x_m} \psi_{d_j(x)}(y) \biggr) \vphi(x-y) \,dy
\biggr\rrvert
\nonumber\\
&&\qquad\leq\biggl\llvert\int_{B(0, d_j(x)/2)} \biggl( \biggl(
\frac{\prt}{\prt x_i} \,\frac{\prt}{\prt x_m} d_j(x) \biggr)
\psi'_{d_j(x)}(y) \biggr) \vphi(x-y) \,dy \biggr\rrvert
\\
&&\qquad\quad{} +\biggl\llvert\int_{B(0, d_j(x)/2)} \biggl( \biggl(
\frac{\prt}{\prt
x_m} d_j(x) \biggr) \biggl( \frac{\prt}{\prt x_i}
d_j(x) \biggr) \psi''_{d_j(x)}(y)
\biggr) \vphi(x-y) \,dy \biggr\rrvert.
\nonumber
\end{eqnarray}
Recall that
the random variables $A^+_x$ and $A^-_x$ are defined
in step 3.
Recall also
that $|\frac{\prt}{\prt x_i} d_j(x)| \leq1$ and
$|\frac{\prt}{\prt x_i} \,\frac{\prt}{\prt x_m} d_j(x)|
\leq c' d_j(x)^{-1}$
for all $i,m$ and $x$. We obtain
\begin{eqnarray*}
&&
\biggl\llvert\int_{B(0, d_j(x)/2)} \biggl( \biggl(\frac{\prt}{\prt x_i}
\,\frac{\prt}{\prt x_m} d_j(x) \biggr) \psi'_{d_j(x)}(y)
\biggr) \vphi(x-y) \,dy \biggr\rrvert
\\
&&\qquad \leq\biggl\llvert\int_{B(0, d_j(x)/2)} \bigl(c'
d_j(x)^{-1} \psi'_{d_j(x)}(y) \bigr)
\vphi(x-y) \,dy \biggr\rrvert
\\
&&\qquad \leq c_5 d_j(x)^{-2} \bigl\llvert\E\bigl(
\vphi\bigl(x-A^+_x\bigr) - \vphi\bigl(x-A^-_x\bigr)
\bigr) \bigr\rrvert.
\end{eqnarray*}
The same reasoning as in (\ref{a146}) and (\ref{a143}) yields
%
\begin{eqnarray}
\label{a144}
&&
\biggl\llvert\int_{B(0, d_j(x)/2)} \biggl( \biggl(
\frac{\prt}{\prt x_i} \,\frac{\prt}{\prt x_m} d_j(x) \biggr)
\psi'_{d_j(x)}(y) \biggr) \vphi(x-y) \,dy \biggr\rrvert
\nonumber\\[-8pt]\\[-8pt]
&&\qquad \leq c_6 d_j(x)^{-d-1} \int
_{B(0, d_j(x)/2)} \bigl|\nabla\vphi(x-z)\bigr| \,dz.
\nonumber
\end{eqnarray}
We apply the same argument with $\psi''$ in place of $\psi'$.
Let $b^+_x( \cdot) = (\psi''_{d_j(x)}( \cdot) \lor0)/\|\psi
''_{d_j(x)}\lor0\|_1$ and
$b^-_x( \cdot) = -(\psi''_{d_j(x)}( \cdot) \land0)/\|\psi
''_{d_j(x)}\land0\|_1$. The functions $b^+_x( \cdot)$ and $b^-_x(
\cdot)$ are probability density functions
that vanish outside the ball $B(0, d_j(x)/2)$.
Let $B^+_x $ and $B^-_x$ be independent
$\RR^d$-valued
random variables with densities $b^+_x( \cdot)$ and $b^-_x( \cdot
)$. We have
\begin{eqnarray*}
&&
\biggl\llvert\int_{B(0, d_j(x)/2)} \biggl( \biggl( \frac{\prt}{\prt x_m}
d_j(x) \biggr) \biggl( \frac{\prt}{\prt x_i} d_j(x) \biggr)
\psi''_{d_j(x)}(y) \biggr) \vphi(x-y) \,dy \biggr
\rrvert
\\
&&\qquad\leq c_7 d_j(x)^{-2} \bigl\llvert\E\bigl(
\vphi\bigl(x-B^+_x\bigr) - \vphi\bigl(x-B^-_x\bigr)
\bigr) \bigr\rrvert
\\
&&\qquad\leq c_7 d_j(x)^{-2}\E\bigl\llvert\vphi
\bigl(x-B^+_x\bigr) - \vphi\bigl(x-B^-_x\bigr) \bigr
\rrvert
\\
&&\qquad\leq c_7 d_j(x)^{-2}\E\int
_{B^+_x}^{B^-_x} \bigl|\nabla\vphi(x-z)\bigr| \,dz,
\end{eqnarray*}
where the last integral is along a line segment from
$B^+_x$ to $B^-_x$.
By step 3, the measure $\E\int_{B^+_x}^{B^-_x} dz$ has a density
that is bounded above by $c_8 d_j(x)^{1-d}$ and vanishes outside the
ball $B(0, d_j(x)/2)$.
In other words, for every set $K\subset D$, $\E\int_{B^+_x}^{B^-_x}
\bone_K(z) \,dz\leq c_8 d_j(x)^{1-d}
m(K\cap B(0, d_j(x)/2))$.
It follows that
%
\begin{eqnarray}
\label{a149}
&&\biggl\llvert\int_{B(0, d_j(x)/2)} \biggl( \biggl(
\frac{\prt}{\prt x_m} d_j(x) \biggr) \biggl( \frac{\prt}{\prt x_i}
d_j(x) \biggr) \psi''_{d_j(x)}(y)
\biggr) \vphi(x-y) \,dy \biggr\rrvert
\nonumber\\
&&\qquad\leq c_7 d_j(x)^{-2}\E\int
_{B^+_x}^{B^-_x} \bigl|\nabla\vphi(x-z)\bigr| \,dz
\nonumber\\[-8pt]\\[-8pt]
&&\qquad\leq c_7 d_j(x)^{-2} \int
_{B(0, d_j(x)/2)} \bigl|\nabla\vphi(x-z)\bigr| c_8
d_j(x)^{1-d} \,dz
\nonumber
\\
&&\qquad= c_9 d_j(x)^{-d-1} \int_{B(0, d_j(x)/2)}
\bigl|\nabla\vphi(x-z)\bigr| \,dz.
\nonumber
\end{eqnarray}
We combine (\ref{a172}), (\ref{a142}), (\ref{a151}), (\ref
{a144}) and (\ref{a149}), and then we use
H\"older's inequality to see that
%
\begin{eqnarray}
\label{a121}
&&
\biggl\llvert\frac{\prt}{\prt x_i} \biggl(\frac{\prt
}{\prt x_m}
\vphi_j(x) \biggr)\biggr\rrvert\nonumber\\
&&\qquad\leq c_{10}
d_j(x)^{-d-1} \int_{B(x, d_j(x)/2)} \bigl|\nabla
\vphi(z)\bigr| \,dz
\\
&&\qquad\leq c_{11} d_j(x)^{-d-1} d_j(x)^{d/2}
\biggl(\int_{B(x, d_j(x)/2)} \bigl|\nabla\vphi(z)\bigr|^2 \,dz
\biggr)^{1/2}
\nonumber
\\
&&\qquad= c_{11} d_j(x)^{-d/2-1} \biggl(\int
_{B(x, d_j(x)/2)} \bigl|\nabla\vphi(z)\bigr|^2 \,dz \biggr)^{1/2}.
\nonumber
\end{eqnarray}

Recall that $\widehat Q = \bigcup_{x\in Q} B(x, d_j(x)/2)$.
We will prove that
%
\begin{equation}
\label{a1421} \int_{\widehat Q} \bigl|\nabla(\vphi)\bigr|^2 (x)
\,dx \geq c_{12} 2^{k(2-d)}a
\end{equation}
for some constant $c_{12}$,
where $a$ is the constant defined in (\ref{e216}).
If the inequality holds with $c_{12}=1$, then we are done. So let us
suppose that
%
\begin{equation}
\label{a1422} \int_{\widehat Q} \bigl|\nabla(\vphi)\bigr|^2 (x)
\,dx \leq2^{k(2-d)}a.
\end{equation}
We combine this with (\ref{a121}) to see that for $x\in Q$,
\begin{eqnarray*}
\biggl\llvert\frac{\prt}{\prt x_i} \biggl(\frac{\prt}{\prt x_m} \vphi_j(x)
\biggr)\biggr\rrvert
&\leq& c_{11} d_j(x)^{-d/2-1}
\biggl(\int_{B(x, d_j(x)/2)} \bigl|\nabla\vphi(z)\bigr|^2 \,dz
\biggr)^{1/2}
\\
&\leq& c_{11} d_j(x)^{-d/2-1} \biggl(\int
_{\widehat Q} \bigl|\nabla\vphi(z)\bigr|^2 \,dz \biggr)^{1/2}
\\
&\leq& c_{11} d_j(x)^{-d/2-1} \bigl(2^{k(2-d)}a
\bigr)^{1/2}.
\end{eqnarray*}
It follows from this and (\ref{eA43}) that the set of $x\in Q_n$
such that $|\nabla\vphi_j (x)| \geq a_n^{1/2}
2^{k}/2$ contains a ball with radius greater than
\[
c_{13} a_n^{1/2} 2^{k} /
\bigl(d_j(x)^{-d/2-1} \bigl(2^{k(2-d)}a
\bigr)^{1/2} \bigr) = c_{13} a_n^{1/2}a^{-1/2}2^{kd/2}d_j(x)^{d/2+1}
\]
and, therefore, it has a volume greater than
\[
\bigl(c_{13} a_n^{1/2}a^{-1/2}
2^{kd/2}d_j(x)^{d/2+1} \bigr)^d =
c_{14} a_n^{d/2}a^{-d/2}2^{kd^2/2}d_j(x)^{d^2/2+d}.\vadjust{\goodbreak}
\]
Hence
\begin{eqnarray*}
\int_{Q_n} \bigl|\nabla\vphi_j (x)\bigr|^2 \,dx
&\geq&\bigl(a_n^{1/2} 2^{k} /2
\bigr)^2 c_{14} a_n^{d/2}a^{-d/2}2^{kd^2/2}d_j(x)^{d^2/2+d}
\\
&=& c_{15} a_n^{1+d/2}a^{-d/2}2^{k(2+d^2/2)}d_j(x)^{d^2/2+d}
\end{eqnarray*}
and, therefore,
%
\begin{eqnarray}
\label{a131} \int_{Q} \bigl|\nabla\vphi_j
(x)\bigr|^2 \,dx &=& \sum_{n=1}^N
\int_{Q_n} \bigl|\nabla\vphi_j(x )\bigr|^2 \,dx
\nonumber\\[-8pt]\\[-8pt]
&\geq&\sum_{n=1}^N c_{15}
a_n^{1+d/2}a^{-d/2}2^{k(2+d^2/2)}d_j(x)^{d^2/2+d}.
\nonumber
\end{eqnarray}
By the H\"older inequality and (\ref{a125}),
\begin{eqnarray*}
\sum_{n=1}^N a_n^{1+d/2}
&\geq& N^{-d/2} \Biggl(\sum_{n=1}^N
a_n \Biggr)^{1+d/2} \geq N^{-d/2} \bigl(a
d^{-1} 2^{1-d} \bigr)^{1+d/2} \\
&=& c_{16}
N^{-d/2} a^{1+d/2}.
\end{eqnarray*}
This and (\ref{a131}) give
\begin{eqnarray*}
\int_{Q} \bigl|\nabla\vphi_j(x)\bigr|^2
\,dx &\geq& c_{17} a N^{-d/2} 2^{k(2+d^2/2)}d_j(x)^{d^2/2+d}\\
&=& c_{17} a \bigl(2^{(k-\ell)d} \bigr)^{-d/2}
2^{k(2+d^2/2)}d_j(x)^{d^2/2+d}
\\
& = & c_{17} a 2^{2k} d_j(x)^{d}
\bigl(d_j(x)2^\ell\bigr)^{d^2/2} \\
&\geq&
c_{18} a 2^{2k} d_j(x)^{d} \geq
c_{19}a 2^{k(2 -d)}
\\
&=& c_{19} 2^{k(2-d)} \sum_{x,y \in(2^{-k}\Z)^d \cap Q, |x- y|=2^{-k}}
\bigl(\vphi_j(x) - \vphi_j(y)\bigr)^2.
\end{eqnarray*}
It follows from this and (\ref{a1420}) that
\begin{eqnarray*}
2^{k(2-d)} \sum_{x,y \in(2^{-k}\Z)^d \cap Q, |x- y|=2^{-k}} \bigl(
\vphi_j(x)- \vphi_j(y)\bigr)^2
&\leq&
c_{20} \int_{Q} \bigl|\nabla\vphi_j(x)\bigr|^2 \,dx \\
&\leq& c_{21} \int_{\widehat Q} \bigl|\nabla
\vphi(x) \bigr|^2 \,dx.
\end{eqnarray*}
In view of (\ref{a1421}) and (\ref{a1422}), we conclude that the last
inequality is always valid.
Recall the constant $M$ from step 5.
Summing over all $Q \in\bigcup_{\ell\leq k}\M_\ell$, $Q \subset
D_R^c$, we obtain
\[
2^{k(2-d)} \sum_{\overline{xy} \in\A'_k, \overline{xy} \subset D_R^c}
\bigl(
\vphi_j(x) - \vphi_j(y)\bigr)^2 \leq M
c_{21} \int_{D_R^c} \bigl|\nabla\vphi(x)\bigr|^2
\,dx \leq M c_{21} \eps.
\]
This shows that (\ref{eA83}) holds and completes the proof of the
theorem.
\end{pf}

\begin{example}\label{a161}
Let $C^1_b(D)$ be the family of bounded continuous functions on $D$
with continuous bounded first order derivatives. Using mean value
theorem, it is easy to see that the inequality (\ref{e21}) holds for
every $\varphi\in C^1_b(D)$ [in fact equality holds for such $\varphi$
since $|\nabla\varphi|$ is bounded on $D$ and so $\lim_{R\to\infty}
\int_{D\setminus D_R}|\nabla\varphi(x)|^2 \,dx=0$]. However, we will
sketch an example, without proof, of a domain $D$ such that $C^1_b(D)$
is not dense in $W^{1,2}(D)$. The point of this example is to show that
Theorem~\ref{T21} cannot be strengthened by adding an extra property
that the functions $\{\vphi_j\}_{j\geq1}$ belong to $C^1_b(D)$.

Let
\begin{eqnarray*}
D_- &=& \bigl\{(x_1,x_2)\in\R^2\dvtx -1 <
x_1 < 0, 0 < x_2 < 1\bigr\},
\\
D_+ &=& \bigl\{(x_1,x_2)\in\R^2\dvtx 0 <
x_1 < 1, 0 < x_2 < 1\bigr\},
\\
D_n &=& \bigl\{(x_1,x_2)\in\R^2\dvtx
-1/n < x_1 < 1/n, 1/n < x_2 < 1/n + \delta_n
\bigr\},
\\
\prt D_n^+ &=& \bigl\{(x_1,x_2)\in
\R^2\dvtx -1/n < x_1 < 1/n, x_2 = 1/n +
\delta_n\bigr\},
\\
\prt D_n^- &=& \bigl\{(x_1,x_2)\in
\R^2\dvtx -1/n < x_1 < 1/n, x_2 = 1/n \bigr\},
\\
D &=& D_- \cup D_+ \cup\bigcup_{n= 2}^\infty
D_n \Bigm\backslash\bigcup_{n= 2}^\infty
\bigl(\prt D_n^+ \cup\prt D_n^-\bigr) .
\end{eqnarray*}
We choose $\delta_n>0$ so small that $D_n$'s are disjoint.
Consider a continuous function $\vphi$ such that $\vphi(x) = 1
$ for $x \in D_+ \setminus\bigcup_{n\geq2} D_n$, $\vphi(x) =
-1 $ for $x \in D_+ \setminus\bigcup_{n\geq2} D_n$ and
$\vphi$ is linear in every $D_n$. The widths $\delta_n$ of
``channels'' $D_n$ can be chosen so small that $\vphi\in
W^{1,2}(D)$ and, moreover, $\int_D |\nabla\vphi|^2 $ can be
made arbitrarily small.

We claim that the function $\vphi$ described above cannot be
approximated by
functions $\eta\in C^1_b(D)$
with arbitrary accuracy. The reason is that for any such~$\eta$,
the oscillation of $\eta$ in a set $D_n$ is arbitrarily small,
for large $n$. Hence, in a neighborhood of $(0,0)$, either
$|\vphi- \eta|$ is nonnegligible on a nonnegligible set,
or $|\nabla\vphi- \nabla\eta|$ is nonnegligible on a
nonnegligible set.
We leave the details to the reader because the claim made in this
example is
not needed for our main theorem.
\end{example}

\section{Invariance principle for reflecting random walk}\label{S3}

Let $\CC$ be the algebra generated by functions
$\{\vphi_j\}_{j\geq1}$ from Theorem~\ref{T21} over $\QQ$. By
the same proof as that for Lemma 2.2 in~\cite{C}, we have the
following.

\begin{lemma}\label{L22} There exists a metric $\rho$ on $D$
which induces
the same Euclidean topology inside $D$
and such that the $\rho$-completion
$\widetilde D$
of $D$ is a regularizing space for Dirichlet form $(\EE, W^{1,2}(D))$.
Moreover, $\CC$ is dense in $C_b(\widetilde D, \| \cdot\|_\infty)$.
\end{lemma}

Let\vspace*{1pt} $m$ be the Lebesgue measure on $D$ extended to
$\widetilde D$ by setting $m(\widetilde D \setminus D)=0$. Then $(\EE,
W^{1,2}(D))$ is a strongly local regular Dirichlet form on
$L^2(\widetilde D; m)$. Let $\widetilde X$ be the Hunt process on
$\widetilde D$ associated with the regular Dirichlet form $(\EE,
W^{1,2}(D))$ on $L^2(\widetilde D; m)$, which is continuous and has
infinite lifetime. Denote by $j$ the projection map from $\widetilde D$
to $\overline D$. Then $X:=j(\widetilde X)$ is a continuous Markov
process taking values on $\overline D$. In general $X$ may not be a
strong Markov process, as one can see from the example when $D$ is the
unit disk in $\RR^2$ with a slit $(-1, 0)\times\{0\}$ removed. Both
$\widetilde X$ and $X$ can be called the reflecting Brownian motion on
$D$.

We will now discuss the relationship between reflecting Brownian
motion $\widetilde X$ on $\widetilde D$ and a better known construction of
reflecting Brownian on an arbitrary domain~$D$.
For an arbitrary
bounded domain $D$ in $\RR^d$,
Fukushima
\cite{Fu} used the Martin--Kuramochi compactification $D^*$ of
$D$ to construct a
conservative continuous Hunt process taking values
in $D^*$.
The process $X^*$ is associated with the regular Dirichlet form
$(\EE, W^{1,2}(D))$
on $L^2( D^*; m^*)$, where $m^*$ is
Lebesgue measure on $D$ extended to $D^*$ by setting
$m^*(D^*\setminus D)=0$.
Since each $f_i(x)\df x_i$ is a function in $W^{1,2}(D)$, it admits
a quasi-continuous extension to $D^*$, which we still denote by $f_i$.
These functions
induce a quasi-continuous projection map $j^*=(f_1, \ldots, f_d)$
from $ D^*$ to $\overline D$.
Then $X':=j^*( X^*)$ is a continuous Markov process taking values on
$\overline D$,
which is called reflecting Brownian motion on $\overline D$ in~\cite{C}.
Both $\widetilde D$ and $D^*$ are regularizing spaces for the
Dirichlet form $(\EE, W^{1,2}(D))$ and so $X$ and $X'$ have
the same finite-dimensional distributions under the initial distribution
$m$. For $x\in D$, both $X$ and $X'$ starting from $x$ behave
like Brownian motion before they hit the boundary after a positive
period of time. Consequently, $X$ and $X'$ have the same finite-dimensional
distributions starting from any interior point
in~$D$. We
can consider processes $X$ and $X'$ as maps from their
underlying probability spaces into the space of continuous
functions $C([0,\infty); \R^d)$. Then the distributions of $X$
and $X'$
in $C([0,\infty); \R^d)$ are identical,
either with initial distribution $m$ or with initial starting point
in $D$. In this sense,
convergence of reflecting random walks to
$X$ or $X'$ is an equally strong result.

Without loss of generality, we assume that $D$ contains the origin $0$.
Recall the definition of $D_k$ from the previous section.
We view $(2^{-k}\Z)^d\cap D_k$ as a graph whose vertices
are $(2^{-k}\Z)^d \cap D_k$, and there is an edge between two vertices
$x$ and $y$ if and only $|x-y|=2^{-k}$ and the line segment
connecting $x$ and $y$ is contained in $D_k$.
By abuse of notation, in this section we will use $D_k$ to denote
the graph $(2^{-k}\Z)^d \cap D_k$.\vadjust{\goodbreak}

For
$x\in D_k$, we use $v_k(x)$ to denote the degree of the vertex $x$
in $D_k$. Let $\{X^k_{j2^{-2k}}, j=0, 1, \ldots\}$ be the simple
random walk
on $D_k$ that jumps every $2^{-2k}$
units of time.
By definition, the random walk $\{X^k_{j2^{-2k}}, j=0, 1, \ldots\}$
jumps to one of its nearest neighbors in $D_k$
with equal probabilities. This
discrete time Markov chain is symmetric with respect to measure
$m_k$, where $m_k(x)=\frac{v_k(x)}{2d} 2^{-kd}$ for $x\in D_k$.
Clearly $m_k$ converge weakly to $m$ on $D$. We now extend the
time-parameter of $ \{X^k_{j2^{-2k}}, j=0, 1, \ldots\}$
to all nonnegative reals using linear interpolation over the
intervals $((j-1)2^{-2k}, j2^{-2k})$ for $j=1, 2, \ldots\,$. We thus
obtain a process
$ X^k=
\{X^k_t, t\geq0\}$. Its law with $X^k_0=x$ will
be denoted by $\P^k_x$.

Recall $\A'_k$ from the beginning of Section~\ref{S2}.
Let $Q_k(x, dy)$ denote the one-step
transition probability for the discrete time Markov chain
$ \{X^k_{j2^{-2k}},\break j=0, 1, \ldots\}$; that is, for
$f\geq0$ on $D$ and $x\in D_k$,
\[
Q_k f(x):=\int_D f(y) Q_k(x,
dy):= \frac1{v_k(x)} \sum_{y\in D_k\dvtx \overline{xy} \in\A'_k}
f(y).
\]

For $f\in C^2(\overline D)$, define
\begin{eqnarray*}
{\mathcal L}_k f(x) :\!&=& \int_D \bigl(f(y)-f(x)
\bigr)Q_k(x, dy)\\[-2pt]
&=& \frac1{v_k(x)} \sum
_{y\in D_k\dvtx \overline{xy} \in\A'_k} \bigl(f(y)-f(x)\bigr),\qquad x\in D_k.
\end{eqnarray*}
Then $ \{f(X^k_{j2^{-2k}}) - \sum_{i=0}^{j-1} {\mathcal L}_k
f(X^k_{i2^{-2k}}), \GG^k_{j2^{-k}}, j=0, 1, \ldots\}$ is a
martingale for every $f\in C^2(\overline D)$, where $\GG^k_t:=\sigma
( X^k_s, s\leq t)$.\vspace*{1pt}

To study the weak limit of $\{X^k, k\geq1\}$, we introduce an
auxiliary process $Y^k$ defined by
$Y^k_t:=X^k_{[2^{2k}t]2^{-2k}}$,
where $[\alpha]$ denotes the largest integer that is less than or
equal to $\alpha$. Note that $Y^k$ is a time-inhomogeneous Markov
process. For every fixed $t>0$, its transition probability\vspace*{1pt}
operator is symmetric with respect to the measure $m_k$ on $D_k$.
Let $\FF^k_t:=\sigma(Y^k_s, s\leq t)$.
By abuse of notation, the
law of $Y^k$ starting from $x\in D_k$ will also be denoted by
$\P^k_x$.

Note that\vspace*{1pt} $Y^k_t = X^k_t$ for every $t$ of the form $t= j2^{-2k}$,
where $j$ is an integer.
Moreover, $\sup_{t\geq0} |X^k_t - Y^k_t| \leq2^{-k}$. It follows that
if the laws of one of the sequences $\{X^k, k\geq1\}$ or $\{Y^k, k\geq
1\}$
converge to a limit on $\D([0, T], \widetilde D)$ for some~$T$,
then the same holds for the other sequence.\vspace*{-3pt}

\begin{theorem}\label{T31}
Let $D$ be a bounded domain in $\RR^n$. Then the laws
$\{\P^k_{m_k},\break
k\geq1\}$ of $\{Y^k, k\geq1\}$ are tight in the space $\D([0, T],
\widetilde D)$ for every $T>0$.\vspace*{-3pt}
\end{theorem}

\begin{pf}
Without loss of generality, we assume that $T=1$. By~\cite{EK},
Theorem~3.9.1, and Lemma~\ref{L22}, it suffices to show that for every
$g\in\CC$, $ \{g(X^{k}) \}_{k \ge1}$ is relatively compact in $\D([0,
1], \RR)$ with the initial distribution $\P^{k}_{m_k}$.\vadjust{\goodbreak}

For each fixed $k\geq1$, we may assume, without loss of
generality, that $\Omega$ is the canonical space $\D([0, \infty),
\widetilde D)$, and $Y^k_t$ is the coordinate map on $\Omega$. Given $t>0$
and a path $\omega\in\Omega$, the time reversal operator $r_t$ is
defined by
%
\begin{equation}
\label{eqntimerevers} r_t(\omega) (s):= \cases{ \omega\bigl((t-s)-
\bigr), &\quad if $0\le s\le t$,
\vspace*{2pt}\cr
\omega(0), &\quad if $s\ge t$.}
\end{equation}
Here for $r>0$, $\omega(r-):=\lim_{s\uparrow r} \omega(s)$ is the
left limit at $r$, and we use the convention that $\omega
(0-):=\omega(0)$. We note that
%
\begin{equation}
\label{eqntimerevleft}
\lim_{s\downarrow0} r_t (\omega) (s)=
\omega(t-) =r_t(\omega) (0) \quad\mbox{and}\quad
\lim_{s\uparrow t}
r_t(\omega) (s) =\omega(0)=r_t (\omega) (t).\hspace*{-28pt}
\end{equation}
Observe that for every integer $T\geq1$, $\P^{k}_{m_k}$ restricted to
the time interval $[0, T)$ is invariant under the time-reversal
operator $r_T$. Note that
\[
M^{k,f}_t:=f\bigl(Y^k_t\bigr) - f
\bigl(Y^k_0\bigr) - \sum_{i=0}^{[2^{2k}t]-1}
{\mathcal L}_k f\bigl(Y^k_{i2^{-2k}}\bigr)
\]
is an $\{\FF^k_t, t\geq0\}$-martingale for every $f\in
\CC$ (cf.~\cite{CFKZ}).
We have
%
\begin{equation}
\label{eqna1} f\bigl(Y^k_t\bigr)-f\bigl(Y^k_0
\bigr)=\tfrac12 M^{k,f}_t -\tfrac12 \bigl( M^{k,f}_{T-}-M^{k, f}_{(T-t)-}
\bigr) \circ r_{T} \qquad\mbox{for } t\in[0, T).\hspace*{-32pt}
\end{equation}

For each $M^{k,f}$, there
exists a continuous predictable
quadratic variation process $\langle M^{k,f} \rangle_t$.
Note that (e.g., see page 214 of~\cite{FOT})
\begin{eqnarray*}
\bigl\langle M^{k,f} \bigr\rangle_t-\bigl\langle
M^{k,f} \bigr\rangle_s &=& \int_s^t
\sum_{y \in D_k} \bigl(f\bigl(Y^{k}_u
\bigr)-f(y)\bigr)^2 Q_k\bigl(Y^{k}_u,
y\bigr)m_k(y)\,du
\\
&\leq& 2d \|f\|_\infty^2 (t-s).
\end{eqnarray*}
Thus by Proposition VI.3.26 in~\cite{JS}, $\{\langle M^{k,f} \rangle_t\}
_{k \ge1}$ is $\CC$-tight
in $\D([0, 1], \RR)$.
As $m_k$ converges weakly to $m$ on $\widetilde D$,
by~\cite{JS}, Theorem VI.4.13, the laws of
$\{ M^{k,f} \}_{k \ge1}$ are
tight in the sense of Skorokhod topology on
$\D([0, 1], \RR)$
with the initial distribution $\P^{k}_{m_k}$.
Since
the laws of $\{ M^{k,f}, t\in[0, 1], \P^{k}_{m_k} \}_{k \ge1}$
are the same as the laws of $\{ M^{k,f}_{(1-t)}, t\in[0, 1], \P
^{k}_{m_k} \}_{k \ge1}$,
it follows from (\ref{eqna1}) that $ \{f(X^{k}) \}_{k \ge1}$
and, consequently,
$ \{g(X^{k}) \}_{k \ge1}$
is tight (and so relatively compact)
in\vspace*{1pt} the sense of Skorokhod topology on
$\D( [0, 1], \RR)$
with the initial distribution $\P^{k}_{m_k}$.
\end{pf}

Let $(\widetilde X, \P)$ be a subsequential limit
of $\{(Y^k, \P^k_{m_k}); k \geq1\}$ on
$\D([0, T], \widetilde D)$.

\begin{theorem}\label{T32}
$(\widetilde X, \P)$ is a stationary symmetric Markov process.
\end{theorem}

\begin{pf}
Let $(\widetilde X, \P)$ be a subsequential limit
of $\{(Y^{k}, \P^k_{m_k}); k \geq1\}$
on $\D([0, T]$, $\widetilde D)$,
say, along a subsequence $\{n_k, k \geq1\}$.
It suffices to show that the finite-dimensional
distributions of $(\widetilde X, P)$\vadjust{\goodbreak}
are determined by a semigroup.
Clearly, $\widetilde m$
is an invariant measure for $\widetilde X$.
For every $t\in[0, T]$, define a linear bounded operator on
$L^2(D; m)=L^2(\widetilde D; \widetilde m)$ by
\[
P_t f \df\E\bigl[ f(\widetilde X_t) \mid\widetilde
X_0 \bigr],\qquad f\in L^2(D; m).
\]
We are going to show that $\{ P_t, t\geq0\}$ is a strongly
continuous symmetric semigroup on $L^2(D; m)$.

(i) We first show that each $P_t$ is a bounded symmetric operator
on $L^2(D; m)$. For every
$f, g\in C_b(\widetilde D, \rho)$ and $t>0$,
it follows from the symmetry of $ (Y^k, \P^k_{m_k})$ that
%
\begin{eqnarray}\label{e34}\quad
\int_D f(x) P_t g (x) m(dx) &=& \E\bigl[ f(
\widetilde X_0) g(\widetilde X_t) \bigr] =
\lim_{k\to\infty} \E_{m_{n_k}} \bigl[ f\bigl(Y^{n_k}_0
\bigr) g\bigl( Y^{n_k}_t\bigr) \bigr]
\nonumber\\
&=&\lim_{k\to\infty} \E_{m_{n_k}} \bigl[ g\bigl(Y^{n_k}_0
\bigr) f\bigl( Y^{n_k}_t\bigr) \bigr] = \E\bigl[ g(
\widetilde X_0) f(\widetilde X_t) \bigr] \\
&=&\int
_D g(x) P_t f (x) m(dx).\nonumber
\end{eqnarray}
In particular, by taking $g=1$, we have
%
\begin{equation}
\label{e35} \int_D \P_t f (x) m(dx) = \int
_D f(x) m(dx) \qquad\mbox{for } f\in C_b(\widetilde
D, \rho).
\end{equation}
Note that
$C_c(D)\subset C_b(\widetilde D, \rho)\subset L^2(D; m)$
and $C_c(D)$ is dense in $L^2(D; m)$.
Hence (\ref{e35}) holds for every $f\in L^2(D; m)$.
Consequently, by the definition of $\P_t$ and Jensen's inequality
for conditional expectation,
%
\begin{equation}
\label{e36} \int_D \bigl(P_t f (x)
\bigr)^2 m(dx) \leq\int_D P_t
\bigl(f^2\bigr) (x) m(dx) = \int_D
f(x)^2 m(dx).
\end{equation}
Hence (\ref{e34}) holds for every
$f, g\in L^2(D, m)$;
in other words, for each $t>0$, $P_t$ is a symmetric contraction
operator in $L^2(D; m)$.

(ii) Next we show that $\{P_t, t\geq0\}$ is a semigroup on $L^2(D; m)$.
For $x=(x_1, \ldots, x_d) \in D_k$, let
$U_k(x):=\prod_{i=1}^d [x_i-2^{-k-1}, x_i+2^{-k-1})$
be the half-closed, half-open cube centered at $x$.
We define
an extension operator $E_k\dvtx\break L^2(D_k,m_k) \to L^2(D,m)$ as follows: for
$g\in L^2(D_k, m_k)$,
%
\begin{equation}
\label{e36b} E_kg(z):= \cases{g(x), &\quad for $z\in U_k(x)$
with $x\in D_k$,
\cr
0, &\quad elsewhere.}
\end{equation}
For $f, g\in C_c(D)$ and $t$ of the form $j 2^{-2l}$, by the uniform continuity,
\begin{eqnarray*}
\lim_{k\to\infty} \int_D f(x) E_{n_k}
P^{n_k}_t g(x)m(dx) &=&\lim_{k\to\infty}
\E_{m_k} \bigl[ f\bigl(Y^k_0\bigr) g
\bigl(Y^k_t\bigr) \bigr]=\E_m \bigl[ f(
\widetilde X_0)g(\widetilde X_t)\bigr]
\\
&=& \int_D f(x) P_t g(x) m(dx).
\end{eqnarray*}
Note that
%
\begin{eqnarray}
\label{e37a} \int_D \bigl(E_{n_k}
P^{n_k}_t g(x)\bigr)^2 m(dx) &=& \int
_{D_{n_k}} \bigl(P^{n_k}_t g(x)
\bigr)^2m_{n_k}(dx) \nonumber\\[-8pt]\\[-8pt]
&\leq&\int_{D_{n_k}}
g(x)^2 m_{n_k}(dx)\nonumber
\end{eqnarray}
and that $\lim_{k\to\infty} \int_{D_{n_k}} g(x)^2 m_{n_k}(dx)
=\int_D g(x)^2 m(dx)$ for $g\in C_c(D)$.
Since every $f\in L^2(D; dx)$ can be approximated in $L^2$-norm
by a sequence
$\{f_k, k\geq1\}\subset C_c(D)$,
we deduce from the last three displays and the
Cauchy--Schwarz inequality that
%
\begin{equation}
\label{e37} \lim_{k\to\infty} \int_D f(x)
E_{n_k} P^{n_k}_t g(x)m(dx) = \int
_D f(x) P_t g(x) m(dx)
\end{equation}
for every $f\in L^2(D; m)$ and $g\in C_c(D)$.

We claim that for every $t=j2^{-2l}$,
%
\begin{equation}
\label{e38}
\lim_{k\to\infty} \int_D \bigl|
E_{n_k} P^{n_k}_t g (x) - P_t g
(x)\bigr|^2 m(dx)=0 \qquad\mbox{for every } g\in C_c(D).\hspace*{-32pt}
\end{equation}
For any fixed $x\in D$, there is $r>0$ so that $B(x, 2r)\subset D$.
When $k$ is large enough,
there is a unique $y_k\in D_k$
so that $x\in U_k(y_k)$. We denote this $y_k$ by $\pi_k(x)$.
Let $q^k(t, x, y)$ denote the transition density
with respect to $m_k$ of simple random walk
on $D_k$ killed upon leaving
$B(x, r)$. It follows from
Donsker's invariance
principle and the uniform H\"older continuity (\cite{D}, Proposition 4.1)
for the parabolic functions of the simple random walk on
$2^{-k}\Z^d$
that $q^k (t, \pi_k(x), \pi_k (y))$
converge locally uniformly
in $y\in B(x_0, r)$ to the transition density $q(t, x, y)$ of
Brownian motion $X$ on $\R^d$ with variance $1/(2d)$ killed upon leaving
$B(x, r)$.
For every $\eps>0$, there is $s>0$
so that $\int_{B(x, r)} q(s, x, y) \,dy >1-(\eps/2)$.
Hence for $k$ sufficiently large, we have
\[
\P_{\pi_k(x)} \bigl(Y^{k}_s\in dy\bigr)= q(s, x, y)
m_k(dy) + \mu_k (dy),
\]
where $\mu_k$ is a signed measure with $| \mu_k | (D_k) < \eps$.
It follows from this and (\ref{e37}) that for $g\in C_c(D)$,
\begin{eqnarray*}
&&\limsup_{k\to\infty} \biggl\llvert P^{n_k}_t g\bigl(
\pi_{n_k}(x)\bigr) -\int_D q(s, x, y)
P_{t-s} g(y) \,dy \biggr\rrvert
\\
&&\qquad= \limsup_{k\to\infty} \biggl\llvert P^{n_k}_s
\bigl(P^{n_k}_{t-s} g \bigr) \bigl(\pi_{n_k}(x)\bigr)
-\int_D q(s, x, y)P_{t-s} g(y) \,dy \biggr\rrvert
\\
&&\qquad\leq \limsup_{k\to\infty} \biggl\llvert\int_{D_k} q (s,
x, y) P^{n_k}_{t-s} g(y) m_k(dy) -\int
_D q(s, x, y) P_{t-s} g(y) \,dy \biggr\rrvert\\
&&\qquad\quad{}+
\eps\| g\|_\infty
\\
&&\qquad= \limsup_{k\to\infty} \biggl\llvert\int_{D} q (s, x,
y) E_{n_k} P^{n_k}_{t-s} g(y) m(dy) -\int
_D q(s, x, y) P_{t-s} g(y) \,dy \biggr\rrvert\\
&&\qquad\quad{}+
\eps\| g\|_\infty
\\
&&\qquad= \eps\| g\|_\infty.
\end{eqnarray*}
Since $\eps>0$ is arbitrary, the above yields that
$\{ P^{n_k}_t g( \pi_{n_k}(x)); k\geq1\}$ is a Cauchy sequence,
and so it converges to some value $u(x)$.
This convergence holds for every $x\in D$, so
we have by (\ref{e37}) that
$u = P_t g$
a.e.; that is, $E_{n_k} P^{n_k}_t g(x)= P^{n_k}_t g( \pi_{n_k}(x))$ converges
to $P_t g (x)$ for a.e. $x\in D$.
Hence by the bounded convergence theorem, (\ref{e38}) holds for every
$g\in C_c(D)$.

Abusing the notation a little bit, for $f\in L^2(D; m)$,
we define $\pi_k f$ as
a function in $L^2(D_k, m_k)$ by
%
\begin{equation}
\label{e311} \pi_k f(x)= \frac{ \int_{U_k (x)} f(y) m(dy)}{
m(U_k(x))} \qquad\mbox{for } x\in
D_k.
\end{equation}
Clearly $\pi_k \circ E_k $ is an identity map on $L^2(D_k, m_k)$ and
\[
\int_{D_k} \bigl|\pi_k f(x)\bigr|^2
m_k(dx) \leq\int_D \bigl|f(x)\bigr|^2 m(dx).
\]
Since $C_c(D)$ is dense in $L^2(D; m)$, we have from (\ref{e38})
that
%
\begin{eqnarray}
\lim_{k\to\infty} \int_D \bigl| E_{n_k}
P^{n_k}_t \pi_{n_k} g (x) - P_t g
(x)\bigr|^2 m(dx)=0\nonumber\\[-8pt]\\[-8pt]
&&\eqntext{\mbox{for every } g\in L^2(D; m).}
\end{eqnarray}
It follows then for $g\in L^2(D; m)$,
\[
P_{t+s} g = \lim_{k\to\infty} E_{n_k} P^{n_k}_t
P^{n_k}_s \pi_{n_k} g = \lim_{k\to\infty}
\bigl(E_{n_k} P^{n_k}_t \pi_{n_k}\bigr)
\bigl(E_{n_k}P^{n_k}_s \pi_{n_k}\bigr) g
= P_t P_s g.
\]
This establishes the semigroup property of $\{P_t, t\geq0\}$.

We have now established that $X$ is a stationary symmetric Markov process.
\end{pf}

The following result is needed in the proof of Theorem~\ref{T34}.

\begin{lemma}\label{L32}
In the above setting, for every $f\in C^\infty_c(D)$, the process
$M^f_t:=f(X_t)- f(X_0)-\frac1{2d} \int_0^t \Delta f(X_s) \,ds$ is a
$\P$-square integrable martingale. This in particular implies that
$\{X_t, t<\tau_D, \P\}$ is a Brownian motion killed upon leaving~$D$,
with initial distribution $m_D$ and infinitesimal generator
$\frac1{2d}\Delta$.
\end{lemma}

\begin{pf}
The proof is the same as that for~\cite{BC}, Lemma 2.2.
\end{pf}

The following is Lemma 2.3 of~\cite{BC}.

\begin{lemma}\label{L33}
Let $D$ be a bounded domain in
$\R^d$
and fix $k\geq1$. Then for
every $j\geq1$ and $f\in L^2(D, m_k)$,
\[
\bigl(f-Q_k^{2j} f, f\bigr)_{L^2(D, m_k)} \leq j
\bigl(f-Q_k^2f, f\bigr)_{L^2(D,
m_k)}
\leq2j(f-Q_k f, f)_{L^2(D, m_k)}.
\]
\end{lemma}

We will say that ``$Z_t$ is a Brownian motion running at speed
$1/n$'' if $Z_{nt}$ is the standard Brownian motion, and we will
apply the same phrase to other related process.

By an argument similar to that in~\cite{BC}, Section 2, but with
Theorem~\ref{T31} in place of~\cite{BC}, Lemma 2.2,
and using
Theorem~\ref{T21} in the energy form argument
in the proof of~\cite{BC}, Theorem 2.4, we can establish the following
theorem.

\begin{theorem}\label{T34}
Let $D$ be a bounded domain in $\R^d$
with $m(\partial D)=0$.
Then for every $T>0$, the
laws of $\{X^k, \P^k_{m_k}\}$ converge weakly\vspace*{1pt} in $C([0, T],
\widetilde D)$ to a stationary reflecting Brownian motion on
$\widetilde D$
running at speed $1/d$ whose initial distribution is the Lebesgue measure
in $D$.
Consequently, for every $T>0$, the
laws of $\{X^k, \P^k_{m_k}\}$ converge weakly in $C([0, T],
\overline D)$ to a stationary reflecting Brownian motion on $\overline D$
running at speed $1/d$ whose initial distribution is the Lebesgue measure
in $D$.
\end{theorem}

\begin{pf}
Fix $T>0$. We know from Theorem~\ref{T31} that the laws of $(X^k,
\P^k_{m_k})$ are tight in the space $\D([0, T], \widetilde D)$. Let
$(\widetilde X, \P)$ be any of subsequential limits, say, along
$(X^{k_j}, \P^{k_j}_{m_{k_j}})$. By\vspace*{1pt} Theorem~\ref{T32} and
its proof, $\widetilde X$ is a time-homogeneous Markov process on
$\widetilde D$ with transition\vspace*{1pt} semigroup $\{P_t, t\geq0\}$ that is
symmetric in $L^2(\widetilde D, m)$. Let $\{P^k_t, t\in2^{-k}\ZZ_+\}$
be defined\vspace*{1pt} by $P^k_t f(x):=\E^k_x[ f(X^k_t)]$. For dyadic $t>0$, say,
$t=j_0/2^{2k_0}$ and $f\in\CC$, we have by Theorem~\ref{T21} and Lemma
\ref{L33},
\begin{eqnarray*}
&&
\frac1{t} (f-P_t f, f)_{L^2(\widetilde D; m)} \\
&&\qquad= \frac1{t}
\lim_{j\to
\infty} \bigl(f-P^{k_j}_t f, f
\bigr)_{L^2(D; m_{k_j})}
\\
&&\qquad= \frac{2^{2k_0}}{j_0} \lim_{j\to\infty} \bigl(f-Q_{k_j}^{j_02^{2k_j-2k_0}}
f, f\bigr)_{L^2(D, m_{k_j})}
\\
&&\qquad\leq \limsup_{j\to\infty} \frac{2^{2k_0}}{j_0} j_02^{2k_j-2k_0}
(f-Q_{k_j} f, f)_{L^2(D, m_{k_j})}
\\
&&\qquad= \limsup_{j\to\infty} 2^{(2-d)k_j} \frac1{2d} \sum
_{x\in
D_{k_j}} \sum_{y\in D_{k_j}\dvtx
\overline{xy} \in\A'_{k_j}} \bigl(
f(x)^2- f(x) f(y) \bigr)
\\
&&\qquad= \frac1{4d} \limsup_{j\to\infty} 2^{(2-d)k_j} \sum
_{x, y\in D_{k_j}\dvtx \overline{xy} \in\A'_{k_j}} \bigl(f(x)-f(y)\bigr)^2
\\
&&\qquad\leq \frac1{2d} \int_D \bigl|\nabla f(x)\bigr|^2 \,dx.
\end{eqnarray*}
Let $(\EE, \FF)$ be the Dirichlet form of $\widetilde X$, or equivalently,
of semigroup $\{P_t,\break t\geq0\}$. That is,
\begin{eqnarray*}
\FF&=& \biggl\{f\in L^2(\widetilde D; m )\dvtx\\
&&\hspace*{6.3pt}\sup_{t>0}
\frac1t (f-P_tf, f)_{L^2(\widetilde D, m)}
=\lim_{t\to0} \frac1t
(f-P_tf, f)_{L^2(\widetilde D, m)}<\infty\biggr\},
\\
\EE(f, f)&=& \sup_{t>0} \frac1t (f-P_tf,
f)_{L^2(\widetilde D; m)}\\
&=&\lim_{t\to0} \frac1t (f-P_tf,
f)_{L^2(\widetilde D; m)} \qquad\mbox{for } f\in\FF.
\end{eqnarray*}
Then for $f\in\CC$,
\[
\EE(f, f) =\sup_{t> 0} \frac1{t} (f-P_t f,
f)_{L^2(\widetilde D; m)}\leq\frac1{2d} \int_D \bigl|\nabla
f(x)\bigr|^2 \,dx.
\]
This shows that $f\in\FF$. As $\CC$ is dense in $(W^{1,2}(D), \|
\cdot\|_{1,2})$ in view of Theorem~\ref{T21}, we have
$W^{1,2}(D)\subset\FF$ and
\[
\EE(f, f) \leq\frac1{2d} \int_D \bigl|\nabla
f(x)\bigr|^2 \,dx \qquad\mbox{for every } f\in W^{1,2}(D).
\]
This, Lemma~\ref{L32}
and~\cite{BC}, Theorem 1.1 (or~\cite{CF}, Theorem 6.6.9) imply
that $\FF=W^{1,2}(D)$ and
\[
\EE(f, f) = \frac1{2d} \int_D \bigl| \nabla
f(x)\bigr|^2 \,dx \qquad\mbox{for } f\in W^{1,2}(D).
\]
We deduce then that $X$ is a stationary reflecting Brownian motion on
$D$ running at speed $1/d$.
This proves that $X^k$ converge weakly on
$C([0, T], \widetilde D)$ to the stationary reflecting Brownian motion
on $\widetilde D$
running at speed $1/d$.

The last assertion comes from the fact that the projection map from
$(\widetilde D, \rho)\to\overline D $ is continuous.
\end{pf}

\begin{remark}
Note that for every $x\in D$, $\widetilde X$ starting from $x$
is a Brownian motion in $D$ before hitting the boundary of $D$,
and the mass of $X$ spreads immediately across the whole set
$D$ immediately after the clock starts, while
$X^k$ starting from $x$ runs like a simple random walk on $2^{-k}\ZZ$
before hitting the boundary of $D_k$. Using these properties,
the weak convergence in Theorem~\ref{T34} can be strengthened\vadjust{\goodbreak}
to show that $(X^k, \P^k_x)$ converges weakly to $(\widetilde X, \P_x)$
for every interior starting point $x\in D$.
We leave the details to the reader.
\end{remark}

\section{Continuous-time reflected random walk} \label{S4}

In this section, we show that reflected Brownian motion on $\overline D$
can be approximated by continuous-time random walks on grids.

Let $D$ be a bounded domain in $\RR^d$ and
$D_k$ be defined as in the beginning of Section~\ref{S2}.
But in this section, $X^k$ will be the continuous time simple random
walk on
$D_k$, making jumps at the rate $2^{-2k}$. By definition, $X^k$
jumps to one of its nearest neighbors with equal probabilities. This
process is symmetric with respect to measure $m_k$, where $m_k(x)=
\frac{v_k(x)}{2d}2^{-kd}$
for $x\in D_k$. Note that $m_k$ converge
weakly to the Lebesgue measure $m$ on $D$, and recall
$\A'_k$ from the beginning of Section~\ref{S2}.
The Dirichlet form of $X^k$ on $L^2(D_k; m_k)$
is given by\looseness=-1
%
\begin{equation}
\label{e41} \EE^k(f, f)= \frac1{4d} \sum
_{x, y \in D_k\dvtx\overline{xy} \in\A'_{k}} 2^{-(d-2)k} \bigl
(f(x)-f(y)\bigr)^2.
\end{equation}\looseness=0

Let $\P^k_{m_k}$ denote the distribution of $\{X^k_t, t\geq0\}$
with the initial distribution~$m_k$.

\begin{lemma}\label{L41}
Assume that $D$ is a bounded domain in $\RR^d$.
For every \mbox{$T>0$}, the laws of stationary random walks
$\{X^k, \P^k_{m_k}, k \geq1\}$ are tight in the space $\D([0, T],
\overline D)$ equipped with the Skorokhod topology.
\end{lemma}

\begin{pf}
The proof is the same as that for~\cite{BC}, Lemma 3.2,
so we omit~it.%
\end{pf}

\begin{theorem}\label{T42}
Let $D$ be a bounded domain in $\RR^d$. Then for every $T>0$, the
stationary random walks $X^k$ on $D_k$ converge weakly in the
space $\D([0, T], \overline D)$, as $k\to\infty$, to the
stationary reflected Brownian motion on $\overline D$ running at
speed $1/d$, whose initial distribution is the Lebesgue measure in~$D$.
\end{theorem}

\begin{pf}
Let $(Z, \P)$ be any of the subsequential limits of $(X^k,
\P^k_{m_k})$ in the space $\D([0, T], \widetilde D)$, say, along $X^{k_j}$.
A similar argument as that for Theorem~\ref{T32} shows that $(Z, \P
)$ is a
time-homogeneous Markov process and its transition
semigroup $\{P_t, t\geq0\}$ is symmetric with respect to the
measure $m$ on $D$.
Furthermore, by a similar argument as that in the
proof of~\cite{BC}, Lemma~2.2, the process $Z$ killed upon leaving
$D$ is a killed Brownian motion in $D$ with speed $1/d$. Let
$(\EE, \FF)$ be the Dirichlet form associated with~$Z$, and let
$\{P^k_t, t\geq0\}$ be the transition semigroup for $X^k$. As
$X^{k_j}$ converge weakly to $Z$ in $\D([0, T], \overline D)$,
we have for every $t>0$ and every $f, g\in C (\overline D)$,
\[
\lim_{j\to\infty} \bigl(f, P^{k_j}_t g
\bigr)_{L^2(D_k; m_k)} = (f, P_t g)_{L^2(D; m)}.
\]
Let $E_k\dvtx L^2(D_k; m_k)\to L^2(D; m)$ and $\pi_k\dvtx L^2(D; m)
\to L^2(D_k, m)$ be the extension operator and restriction
operator defined by (\ref{e36b}) and (\ref{e311}), respectively.
Then the last display can be restated as
%
\begin{equation}
\label{e42} \quad\lim_{j\to\infty} \bigl( f, E_{k_j}
P^{k_j}_t \pi_{k_j}g\bigr)_{L^2(D; m)} = (f,
P_t g)_{L^2(D; m)} \qquad\mbox{for } f, g\in C(\overline D).
\end{equation}
Note that
\[
\| E_k f\|_{L^2(D; m)}=\|f\|_{L^2(D_k; m_k)}
\qquad\mbox{for } f\in L^2(D_k; m_k)
\]
and
\[
\| \pi_k g\|_{L^2(D_k; m_k)}\leq\| g\|_{L^2(D; m)}\qquad \mbox{for } g\in
L^2(D; m).
\]
Since $C(\overline D)$ is dense in $L^2(D; m)$ and $P^k_t$ and $P_t$
are contraction operators on $L^2(D_k; m_k)$ and $L^2(D; m)$, respectively,
we deduce from (\ref{e42}) that
%
\begin{eqnarray}
\label{e43} \lim_{j\to\infty} \bigl( f, E_{k_j}
P^{k_j}_t \pi_{k_j}g\bigr)_{L^2(D; m)} = (f,
P_t g)_{L^2(D; m)} \nonumber\\[-8pt]\\[-8pt]
&&\eqntext{\mbox{for every } f, g\in L^2(D;
m).}
\end{eqnarray}
Recall from Theorem~\ref{T21} that $\CC$ is the algebra generated
by functions
$\{\vphi_j\}_{j\geq1}$ over~$\QQ$.
As a special case of
(\ref{e43}), we obtain
for every $t>0$ and $f\in\CC$,
\begin{eqnarray*}
(f, P_t f)_{L^2(D; m)} &=& \lim_{j\to\infty} \bigl( f,
E_{k_j} P^{k_j}_t \pi_{k_j}f
\bigr)_{L^2(D; m)} \\
&=& \lim_{j\to\infty} \bigl( \pi_{k_j}f,
P^{k_j}_t \pi_{k_j}f\bigr)_{L^2(D_k; m_k)}.
\end{eqnarray*}
Since $\CC\subset C_b(D)$ and $D$ is bounded,
\[
\lim_{k\to\infty} \int_{D_k} \bigl|f(x) - \pi_k
f(x) \bigr|^2 m_k(dx)=0 \qquad\mbox{for } f\in\CC.
\]
Hence we conclude that
%
\begin{equation}
\label{e44} (f, P_t f)_{L^2(D; m)} = \lim_{j\to\infty} \bigl(
f, P^{k_j}_t f\bigr)_{L^2(D_k; m_k)} \qquad\mbox{for } f\in\CC.
\end{equation}
Thus we have for every $t>0$ and $f\in\CC$,
\begin{eqnarray*}
\frac1t (f, f-P_t f)_{L^2(D; m)} &=& \lim_{j\to\infty} \frac1t
\bigl( f, f- P^{k_j}_t f\bigr)_{L^2(D_k; m_k)}
\\
&\leq& \liminf_{j\to\infty} \sup_{s>0} \frac1s \bigl( f, f-
P^{k_j}_s f\bigr)_{L^2(D_k; m_k)}
\\
&=& \liminf_{j\to\infty} \EE^{k_j} ( f, f)
\\
&=& \liminf_{j\to\infty} \frac1{4d} \sum_{x, y \in D_k\dvtx\overline
{xy} \in\A'_{k}}
2^{-(d-2)k} \bigl( f(x)- f(y)\bigr)^2
\\
&\leq& \frac1{2d} \int_D \bigl|\nabla f (x)\bigr|^2
m(dx),
\end{eqnarray*}
where in the last inequality we used Theorem~\ref{T21}.
Thus
\begin{eqnarray*}
\EE(f, f) &=&\sup_{t>0} \frac1t (f-P_tf, f)_{L^2(D; m)}
\\
&\leq&\frac1{2d} \int_D \bigl|\nabla f(x)\bigr|^2 m( dx )
\qquad\mbox{for every } f\in\CC.
\end{eqnarray*}
Since $\CC$ is dense in the
Sobolev space $W^{1,2}(D)$ with respect to norm $\| \cdot
\|_{1,2}$, it
follows that $\FF\supset W^{1,2}(D)$ and
\[
\EE(f, f) \leq\frac1{2d} \int_D \bigl|\nabla f (x)
\bigr|^2 m(dx) \qquad\mbox{for every } f\in W^{1,2}(D).
\]
Define
\[
\EE^0(f, g)=\frac1{2 d} \int_D \nabla f (x)
\cdot\nabla g(x) m(dx) \qquad\mbox{for } f, g\in W^{1,2}(D).
\]
Note that $(\EE^0, W^{1,2}(D))$ is the Dirichlet form for the
reflected Brownian motion on $D$ running at speed $1/d$. On the
other hand, as we have observed at the beginning of this proof,
the process $Z$ killed upon leaving $D$ is a killed Brownian
motion in $D$ with speed $1/d$. Therefore according to
\cite{BC}, Theorem 1.1 (or~\cite{CF}, Theorem 6.6.9),
$(\EE, \FF)=(\EE^0, W^{1,2}(D))$. In other words, we
have shown that every subsequential limit of $X^k$ is reflected
Brownian motion on $D$ with initial distribution being the
Lebesgue measure on $D$ and with speed $1/d$. This shows that
$X^k$ converges weakly on the space $\D([0, \infty), \overline D)$
to the stationary reflected Brownian motion $X$ on $D$ running at
speed $1/d$.
\end{pf}



\printaddresses

\end{document}